\documentclass[12pt]{article}
\usepackage{amsmath}
\usepackage{amsfonts}
\usepackage{amssymb}
\usepackage{amsthm}
\usepackage{graphicx}
\usepackage{bbm}
\usepackage{hyperref}
\def\arxiv#1{\texttt{\href{http://front.math.ucdavis.edu/#1}{arXiv:#1}}}
\def\MR#1{\href{http://www.ams.org/mathscinet-getitem?mr=#1}{MR#1}}
\def\article#1{\href{#1}{Article}}
\def\aldous{{\href{http://www.stat.berkeley.edu/users/aldous/}{David Aldous}}}
\def\miermont{{\href{http://www.dma.ens.fr/~miermont/}{Gr\'egory Miermont}}}
\def\pitman{{\href{http://www.stat.berkeley.edu/users/pitman/}{Jim Pitman}}}
\def\dept_stat{{\href{http://www.stat.berkeley.edu/}{Department of Statistics}}}
\def\dma{{\href{http://www.dma.ens.fr/}{DMA}}}
\def\lpma{{\href{http://www.proba.jussieu.fr/}{LPMA}}}

\newcommand{\sth}{{{{\mbox{${\mbox{\scriptsize\boldmath$\theta$}}$}}}}}

\def\re#1{(\ref{#1})}
\newcommand{\Bexc}{B^{\rm exc}}
\newcommand{\Bbr}{B^{\rm br}}
\newcommand{\sfrac}[2]{{\textstyle\frac{#1}{#2}}}
\newcommand{\bp}{{\mathbf p}}
\newcommand{\bq}{{\mathbf q}}

\newcommand{\bt}{{\bf t}}

\newcommand{\bT}{{\bf T}}

\newcommand{\BB}{{\cal B}}
\newcommand{\xbrt}{X^{{\rm br},\sth}}
\newcommand{\btheta}{{\mbox{\boldmath$\theta$}}}
\newcommand{\bTheta}{{\mbox{\boldmath$\Theta$}}}
\newcommand{\Thfinite}{\bTheta_{{\rm finite}}}

\renewcommand{\sp}{\sigma({\bf p})}

\newcommand{\gip}{G_I^{{\bf p}}}
\newcommand{\Fexc}{F^{{\rm exc},\bp}}
\newcommand{\fexp}{F^{{\rm exc},{\bf p}}}
\newcommand{\R}{\mathbb{R}}

\newcommand{\tree} { {\cal T }}
\newcommand{\TT} { {\cal T }}

\def\build#1_#2^#3{\mathrel{
\mathop{\kern 0pt#1}\limits_{#2}^{#3}}}

\newcommand{\cd}{\build\to_{}^d}
\newcommand{\ed}{\build=_{}^d}

\newcommand{\TTth}{{\cal T}^{\sth}}
\newcommand{\thv}{\sth}

\newcommand{\HT}{{\rm ht}}
\newcommand{\psbr}{\psi_{{\bf p}}^{{\rm breadth}}}
\newcommand{\psde}{\psi_{{\bf p}}^{{\rm depth}}}
\def\cq{$\hfill \square$}

\def\d{{\rm d}}

\def\eps{\varepsilon}

\def\en{\end{equation}}
\def\lb#1{\label{#1}}
\def\ba{\begin{eqnarray*}}
\def\ea{\end{eqnarray*}}

\def\rem{\noindent{\bf Remark. }}

\def\cp{\build\to_{}^{p}}

\newcommand{\ind}{{\mathbbm{1}}}

\newtheorem{prp}{Proposition}
\newtheorem{lmm}{Lemma}
\newtheorem{thm}{Theorem}
\newtheorem{crl}{Corollary}

\begin{document}
\title{The exploration process of inhomogeneous continuum random trees,
and an extension of Jeulin's
local time identity
}
\author{\aldous\thanks{\dept_stat, 
University of California, 367 Evans Hall \# 3860,
Berkeley, CA 94720-3860. Research supported in part by N.S.F. 
Grant DMS-0203062} and 
\miermont\thanks{\dma, Ecole Normale Sup\'erieure, and \lpma,
Universit\'e Paris 6. 45 rue d'Ulm, 75230 Paris Cedex 05, France} 
and \pitman\thanks{\dept_stat, University of California,
Berkeley. Research supported in part by N.S.F. Grant DMS-0071448}
}

\maketitle

\begin{abstract}
We study the inhomogeneous continuum random trees (ICRT) that arise as
weak limits of birthday trees. We give a description of
the {\em exploration process}, a function defined on $[0,1]$ that
encodes the structure of an ICRT, and also of its {\em width process}, determining
the size of layers in order of height.  These processes turn out to
be transformations of bridges with exchangeable increments, which have already
appeared in other ICRT related topics such as stochastic additive coalescence.
The results rely on two different constructions of birthday trees from
processes with exchangeable increments, on weak convergence arguments,
and on general theory on continuum random trees.
\end{abstract}

\vspace{0.1in}

\vspace{0.1in}

\newpage
\section{Introduction}
This paper completes one circle of ideas
(describing the inhomogeneous continuum random tree)
while motivated by another (limits of non-uniform random $\bp$-mappings
which are essentially different from the uniform case limit).
Along the way, a curious extension of
Jeulin's result on total local time for
standard Brownian excursion will be established.

Consider a continuous function $f:[0,1] \to [0,\infty)$ which is
an ``excursion'' in the sense
\[ f(0) = f(1) = 0; \quad f(u) > 0, \ 0<u<1 . \]
Use $f$ to make $[0,1]$ into the pseudo-metric space with distance
\begin{equation}
\lb{def-explore}
 d(u_1,u_2):=
(f(u_1) - \inf_{u_1 \leq u \leq u_2} f(u)) +
(f(u_2) - \inf_{u_1 \leq u \leq u_2} f(u)) , \
u_1 \leq u_2 .
\end{equation}
After taking the quotient by identifying points of $[0,1]$ that are at
$d$-pseudo distance $0$, this
space is a {\em tree} in that between any two points there is a
unique path; it carries a {\em length measure} induced by the distance $d$,
and a {\em mass measure}, with unit total mass, induced from Lebesgue
measure on $[0,1]$.  An object with these properties can be
abstracted as a {\em continuum tree}.
Using a random excursion function yields a {\em continuum random tree} (CRT):
Aldous \cite{me55,me56}.
The construction of a continuum random tree $\TT$ via a random function $f$,
in this context called the {\em exploration process} of $\TT$ (in Le Gall et
al.
~\cite{lglj97,DGall02}, it is instead called {\em height process} while
the term exploration process is used for a related measure-valued process),
is not the
only way of looking at a CRT; there are also\\
(a) constructions via line-breaking schemes\\
(b) descriptions via the spanning subtrees on $k$ random points chosen
according to mass measure\\
(c) descriptions as weak or strong $n \to \infty$ limits of rescaled
$n$-vertex discrete random trees.\\
As discussed in \cite{me55,me56} the fundamental example is the
{\em Brownian} CRT, whose exploration process is twice
standard Brownian excursion
(this was implicit in Le Gall \cite{LeG90}), with line-breaking
construction given in Aldous \cite{me46},
spanning subtree description in Aldous \cite{me56} and Le Gall \cite{LeG93},
and weak limit
(for conditional Galton-Watson trees) behavior in \cite{me55,me56}
(see Marckert and Mokkadem \cite{mm01} for recent review).
A more general model, the {\em inhomogeneous continuum random tree} (ICRT)
$\TTth$, arose in Camarri and Pitman \cite{jpmc97b} as a weak limit in a certain model
($\bp$-trees) of discrete random trees.
The definition and simplest description of $\TTth$ is via a line-breaking construction
based on a Poisson point process in the plane (Aldous and Pitman
\cite{jpmc97b,me87}), which we recall below. 
The spanning subtree description is set out in Aldous and Pitman \cite{me88},
and the main purpose of this paper is to complete the description
of $\TTth$ by determining its exploration process (Theorem \ref{T1}).

\subsection{Statement of results}
The parameter space $\bTheta$ of the ICRT $\TTth$ is defined
\cite{me87}
to consist of 
sequences
$\btheta = (\theta_0,\theta_1,\theta_2,\ldots)$
such that\\
(i) $\theta_0 \geq 0; \quad
\theta_1 \geq \theta_2 \geq \ldots \geq 0$;\\
(ii) $\sum_i \theta_i^2 = 1$;\\
(iii)
if $\sum_{i=1}^\infty \theta_i < \infty$ then $\theta_0 > 0$.\\
We will often consider the {\em finite-length} subspace
$\Thfinite$ of $\bTheta$ for which $\theta_i = 0 \ \forall i>I$,
for some $I \geq 0$, calling $I$ the {\em length} of $\btheta$.
Note that $\btheta \in \Thfinite$ can be specified by specifying
a decreasing sequence $(\theta_1,\ldots,\theta_I)$ for which
$\sum_{i=1}^I \theta_i^2 < 1$; then set
$\theta_0 = \sqrt{1 - \sum_{i\geq 1} \theta_i^2} > 0$.

Let $\{(U_i,V_i),i\geq 1\}$ be a Poisson measure on the
first octant $\{(x,y):0\leq y\leq x\}$, with intensity $\theta_0^2$ per
unit area. For every $i\geq 1$ let also $(\xi_{i,j},j\geq 1)$ be a Poisson 
process on the positive real
line with intensity $\theta_i$ per unit length.
The hypotheses on $\btheta$ entail that the set of
points $\{U_i,i\geq 1,\xi_{i,j},i\geq 1,j\geq
2\}$ is discrete and can be ordered as
$0<\eta_1<\eta_2<\ldots$, we call them {\em cutpoints}. It is
easy to see that $\eta_{k+1}-\eta_k\to0$ as $k\to \infty$. By convention let
$\eta_0=0$. Given a cutpoint $\eta_k$, $k\geq 1$, we associate a
corresponding {\em
joinpoint} $\eta^*_k$ as follows. If the cutpoint is of the form $U_i$, then
$\eta^*_k=V_i$. If it is of the form $\xi_{i,j}$ with $j\geq 2$, we let
$\eta^*_k=\xi_{i,1}$. The hypothesis $\theta_0>0$ or $\sum_{i\geq 1}\theta_i=
\infty$ implies that joinpoints are a.s.\ everywhere dense in $(0,\infty)$.

The tree is then constructed as follows. Start with a branch $[0,\eta_1]$, and
recursively,  given the tree  is constructed  at stage  $J$, add the line segment
$(\eta_J,\eta_{J+1}]$ by branching its left-end to the joinpoint $\eta^*_J$
(notice that $\eta^*_J<\eta_J$ a.s.\ so that the construction is indeed
recursive as $J$ increases). When all the branches are attached
to their respective joinpoints, relabel the joinpoint corresponding to some
$\xi_{j,1}$ as joinpoint $j$, and forget other labels (of the form $\eta_i$ or 
$\eta^*_i$). We obtain a metric tree (possibly with marked vertices
$1,2,\ldots$), whose completion we call $\TT^{\sth}$. 

A heuristic description of the structure of the ICRT goes as follows. 
When $\theta_0=1$ and hence $\theta_i=0$ for $i\geq 1$,
the tree is the Brownian CRT, it has no marked vertex and it is a.s.\ binary, 
meaning that branchpoints have degree $3$. It is the only ICRT for which the 
width process  defined below is continuous,  and for which  no branchpoint has
degree more than $3$. When $\btheta\in \Thfinite$
has length $I\geq 1$, the structure looks like that
of the CRT, with infinitely many branchpoints with degree $3$, but there
exist also exactly $I$ branchpoints with infinite degree
which  we  call  {\em hubs},  and  these  are  precisely the  marked  vertices
$1,2,\ldots,I$ corresponding to the joinpoints
$\xi_{1,1},\ldots,\xi_{I,1}$ associated to the
Poisson processes with intensities $\theta_1,\ldots,\theta_I$ defined above. 
The width process defined below has
$I$ jumps with respective sizes $\theta_1,\ldots,\theta_I$, which
occur at distinct times a.s. These jump-sizes can be interpreted as the
{\em local time} of the different hubs -- see remark following Theorem 
\ref{T2}. When $\btheta\notin\Thfinite$, then the hubs become
everywhere dense on the tree. Whether there exists branchpoints with degree
$3$ or not depends on whether $\theta_0\neq 0$ or $\theta_0=0$. Also,
the tree can become unbounded.

It turns out that the relevant exploration process is closely related
to processes recently studied for slightly different purposes.
The Brownian CRT in \cite{me82}, and then the ICRT in \cite{me87}, were used
by Aldous and Pitman
to construct versions of the {\em standard}, and then the {\em general},
{\em additive coalescent}, and its dual
{\em fragmentation process}, which are Markov processes on the state space
$\Delta$ of sequences $\{(x_1,x_2,\ldots): \ x_i \geq 0, \sum_i x_i = 1\}$.
In \cite{me82,me87} the time-$t$ state ${\bf X}(t)$ is specified
as the vector of masses of tree-components in the forest
obtained by randomly cutting the Brownian CRT or ICRT at some rate depending
on $t$.
Bertoin \cite{bertoin00f} gave the following more direct construction.
Let $(\Bexc_s, 0 \leq s \leq 1)$ be standard Brownian excursion.
For fixed $t \geq 0$ consider the process of height-above-past-minimum of
\[ \Bexc_s - ts, \quad 0 \leq s \leq 1 . \]
Then its vector of excursion lengths is $\Delta$-valued, and this process
(as $t$ varies) can be identified with the standard case of the additive coalescent.
More generally, for $\btheta \in \bTheta$ 
consider the ``bridge'' process
\[ \theta_0 \Bbr_s + \sum_{i=1}^\infty \theta_i (\ind_{\{U_i \leq s\}} - s), \quad 0 \leq s \leq 1  \]
where $(U_i)$ are independent random variables with uniform law on $(0,1)$.
Use the {\em Vervaat transform} -- relocate the space-time
origin to the location of the infimum -- to define an
``excursion'' process
$X^\sth
= (X^{\sth}_s, \ 0 \leq s \leq 1)$
which has positive but not negative jumps.
Bertoin \cite{bertoin01eac} used $X^{\sth}$ to construct
the general additive coalescent,
and Miermont \cite{miermont01} continued the study of fragmentation processes
by this method.
In this paper we use $X^{\sth}$ to construct a {\em continuous} excursion
process $Y^{\sth}$; here is the essential idea.
A jump of $X^{\sth}$ at time $t_J$ defines an interval
$[t_J,T_J]$ where
$T_J:= \inf \{t>t_J: \ X^{\sth}_t = X^{\sth}_{t_J -} \}$.
Over that interval, replace $X^{\sth}_s$ by
$X^{\sth}_{t_J -} + X^{\sth}_s - \inf_{t_J \leq u \leq s} X^{\sth}_u$.
Do this for each jump, and let $Y^\sth$ be the resulting process.
To write it in a more compact way, the formula
\begin{equation}\label{morecompact}
Y^{\sth}_s=m\left\{\inf_{u\leq r\leq s}X^{\sth}_r:0\leq u\leq s\right\}
\end{equation}
holds, where $m$ is Lebesgue's measure on $\R$.
Details are given in section \ref{sec-cXZ}.
We can now state our main result.
\begin{thm}\label{T1}
Suppose $\btheta \in \bTheta$ satisfies $\sum_i\theta_i<\infty$. Then
the exploration process of the ICRT $\TT^{\sth}$ is distributed
as $\sfrac{2}{\theta_0^2}Y^{\sth}$. 
\end{thm}

As will be recalled in Sect.\ \ref{con-ptrees}, the precise meaning of this 
theorem is: let $U_1,U_2,\ldots$ be independent
uniform variables on $[0,1]$, independent of $Y^{\sth}$, and as around 
(\ref{def-explore}), replacing $f$ by $\sfrac{2}{\theta_0^2}Y^{\sth}$, 
endow $[0,1]$ with a pseudo-distance $d$, so that the natural 
quotient gives a tree
 $\TT^{2\theta_0^{-2}Y}$ where $Y=Y^{\sth}$. Then for every $J\in\mathbb{N}$,
the subtree spanned
by the root  (the class of $0$) and the (classes  of) $U_1,\ldots,U_J$ has the
same law as the tree  $\TT^{\sth}_J$ obtained by performing the stick-breaking
construction until the $J$-th step. Since  $(U_i,i\geq 1)$ is  a.s.\ dense in
$[0,1]$   and  by  uniqueness   of  the   metric  completion,   $\TT^{\sth}$  and
$\TT^{2\theta_0^{-2}Y}$ indeed encode the same random topological space. 
We also note that our proofs easily extend to showing that the 
hub with extra label $i$ is associated to the class of $t_i$ or $T_i$, and
this class is  exactly $\{s\in[t_i,T_i]:Y^{\sth}_s=Y^{\sth}_{t_i}\}$. To avoid
heavier notations, we  will not take these extra labels  into account from now
on. 

When $\sum_i\theta_i=\infty$, the exploration process of the 
ICRT, if it exists, can be obtained as a certain weak limit of processes of 
the form $\sfrac{2}{(\theta_0^{n})^2}Y^{\sth_n}$ for 
approximating sequences $\theta^{n}\in\Thfinite$, and in particular, when 
$\theta_0>0$ one guesses that the exploration process of $\TT^{\sth}$ will 
still be $\sfrac{2}{\theta_0^2}Y^{\sth}$, but we will not concentrate on this 
in the present paper. 

\rem
Formula (\ref{morecompact}) is inspired by the work of Duquesne and
Le Gall \cite{DGall02}, in which continuum random trees (``L\'evy trees'') are 
built out of sample paths of L\'evy processes. Our work suggest that there are 
many similarities between ICRTs and L\'evy trees. In fact, L\'evy
trees turn out to be ``mixings'' of ICRTs in an analogous way that L\'evy
bridges are mixing of extremal bridges with exchangeable increments. This will
be pursued elsewhere.

In principle Theorem \ref{T1} should be provable within the
continuous-space context,
but we do not see such a direct proof.
Instead we use weak convergence arguments.
As background, there are many ways of coding discrete trees as walks.
In particular, one can construct a Galton-Watson
tree with offspring distribution $\xi$ in terms of an excursion
of the discrete-time integer-valued random walk with step
distribution $\xi - 1$.
In fact there are different ways to implement
the same construction, which differ according to how one chooses to order
vertices in the tree, and the two common choices are the
{\em depth-first} and the {\em breadth-first} orders.
In section \ref{con-ptrees} we give a construction of a random $n$-vertex
$\bp$-tree, based on using $n$ i.i.d. uniform$(0,1)$ random variables
to define an excursion-type function with drift rate $-1$ and
with $n$ upward jumps, and again there are two ways to implement
the construction depending on choice of vertex order.
These constructions seem similar in spirit to, but not exactly
the same as, those used in the {\em server system} construction in
\cite{bertoin01eac} or the {\em parking process} construction in
Chassaing and Louchard \cite{chasslou99}.
When $\btheta\in\Thfinite$, by
analyzing asymptotics of the (appropriately rescaled) discrete excursion
using
depth-first order, in the asymptotic regime where convergence to the ICRT
holds, we get weak convergence to the process $Y^{\sth}$,
and we show that this discrete excursion
asymptotically agrees with $\theta_0^2/2$ times
the discrete exploration process; we extend this to the case
$\sum_i\theta_i<\infty$ by approximating the tree $\TT^{\sth}$ by the tree
$\TT^{\sth^n}$ associated to the truncated sequence
$(\theta_1,\ldots,\theta_n,0,\ldots)$, and
that is the proof of Theorem \ref{T1}.
It is a curious feature of the convergence of approximating $\bp$-trees
to $\TTth$ that the rescaled discrete approximation process converges to
$\sfrac{2}{\theta_0^2}Y^{\sth}$ for a topology which is weaker than the usual
Skorokhod topology. In the course of proving Theorem \ref{T1}, we will
give sufficient conditions for this stronger convergence to happen.

For any continuum tree with mass measure $\mu$,
we can define
\[ \bar{W}(h) = \mu \{x: \HT(x) \leq h \} , \quad h \geq 0  \]
where the {\em height} $\HT(x)$ of point $x$ is just its distance
to the root.
If $\bar{W}(h) = \int_0^h W(y) \ \d y, \quad h \geq 0$ then
$W(y)$ is the ``width'' or ``height profile'' of the tree
(analogous to the size of a particular generation in
a branching process model).
The time-changed function
$(W(\bar{W}^{-1}(u)),
 0 \leq u \leq 1)$
can be roughly interpreted as the width of the layer of
the tree containing vertex $u$, where vertices are labelled
by $[0,1]$ in breadth-first order.
Parallel to (but simpler than) the proof of Theorem \ref{T1}
sketched above, we show
that excursions coding $\bp$-trees using {\em breadth-first} order
converge to $X^{\sth}$,
and agree asymptotically with the height profile
(sizes of successive generations) of the $\bp$-tree.
In other words
\begin{thm}\label{T2}
Let $\btheta \in \bTheta$.
For the ICRT $\TT^{\sth}$ the width process
$W(y)=W^{\sth}(y)$ exists, and
\[ (W^{\sth}((\bar{W}^{\sth})^{-1}(u)), 0 \leq u \leq 1) \ed
 (X^{\sth}(u), 0 \leq u \leq 1). \]
\end{thm}
Qualitatively, in breadth-first traversal of the ICRT, when we encounter
a hub at some $0<u<1$ we expect the time-changed width function
$W(\bar{W}^{-1}(u))$
to jump by an amount representing a ``local time'' measuring relative numbers
of edges at that hub.  Theorem \ref{T2} shows these jump amounts are
precisely the $\theta$-values of the hubs.

When $\sum_i\theta_i<\infty$, combining Theorems \ref{T1} and \ref{T2} 
gives a result whose statement does not involve trees:
\begin{crl}\label{C1}
Let $\btheta \in \bTheta$ satisfy $\sum_i\theta_i<\infty$.
The process $\sfrac{2}{\theta_0^2}Y^{\sth}$ has an occupation density
$(W^{\sth}(y), 0 \leq y < \infty)$
satisfying
\[ (W^{\sth}((\bar{W}^{\sth})^{-1}(u)), 0 \leq u \leq 1) \ed
 (X^{{\sth}}(u), 0 \leq u \leq 1). \]
\end{crl}
Note that the ``Lamperti-type'' relation between $W^{\sth}$ and $X^{\sth}$ 
is easily inverted as
\begin{equation}
 (X^{\sth}_{L^{-1}(y)},y\geq 0) \ed (W^{\sth}(y),y\geq 0), \label{Lamperti}
\end{equation}
where
\[ L(t):=\int_0^t\frac{\d s}{X^{\sth}_s} \in[0,\infty],
\hskip1cm 0\leq t\leq 1. \]
This provides a generalization of the following result of
Jeulin \cite{jeulin85} (see also Biane-Yor
\cite{by88e}),
which from our viewpoint is the Brownian CRT case where $\theta_0 = 1$.
Let $(l_u, 0 \leq u < \infty)$ be occupation density
for $(\Bexc_s, 0 \leq s \leq 1)$.
Then
\[ (\sfrac{1}{2}l_{u/2}, 0 \leq u < \infty) \ed
(B^{\rm exc}_{L^{-1}(u)}, 0 \leq u < \infty) \]
where
$L(t):= \int_0^t \frac{1}{\Bexc_s} \d s$.
One might not have suspected a  possible generalization of this
identity to jump processes without the interpretation provided
by the ICRT.

%

Theorem \ref{T2} has the following other corollary:

\begin{crl}\label{httree}
For any $\btheta\in\bTheta$, the height $\sup_{v\in\TT^{\sth}}\HT(v)$ 
of the ICRT $\TT^{\sth}$ has the same law as
$$\int_0^1\frac{\d s}{X^{\sth}_s}.$$
\end{crl}

\subsection{Discussion}
As formulated above, the purpose of this paper is to prove
Theorems \ref{T1} and \ref{T2} concerning the ICRT.
But we have further motivation.
As ingredients of the proof, we take a known result
(Proposition \ref{PwcICRT})
on weak convergence of random $\bp$-trees to the ICRT,
and improve it to stronger and more informative versions 
(Propositions \ref{pGZ} and \ref{Pkey}).
The Theorems and these ingredients will be used in a sequel
\cite{meAMP} studying asymptotics of random $\bp$-mappings.
By using Joyal's bijection between mappings and trees,
one can in a sense reduce questions of convergence of
$\bp$-mappings to convergence of random $\bp$-trees.
In particular,
under a {\em uniform asymptotic negligibility}
hypothesis which implies that the exploration process
of $\bp$-trees converges to Brownian excursion,
one can use a 
``continuum Joyal functional''
(which takes Brownian excursion to reflecting Brownian motion)
to show \cite{me102} 
that the exploration process of the random $\bp$-mappings
converges to reflecting Brownian bridge.
The results of the present paper give the limit exploration process
$Y^{\sth}$
for more general sequences of $\bp$-trees, and to deduce convergence
of the associated random $\bp$-mappings
we need to understand how the continuum Joyal functional acts
on $Y^{\sth}$.
This is the subject of the sequel \cite{meAMP}.

\section{Constructing $X^{\theta}$ and $Y^{\theta}$}
\label{sec-cXZ}

\begin{figure}
\begin{center}
\input{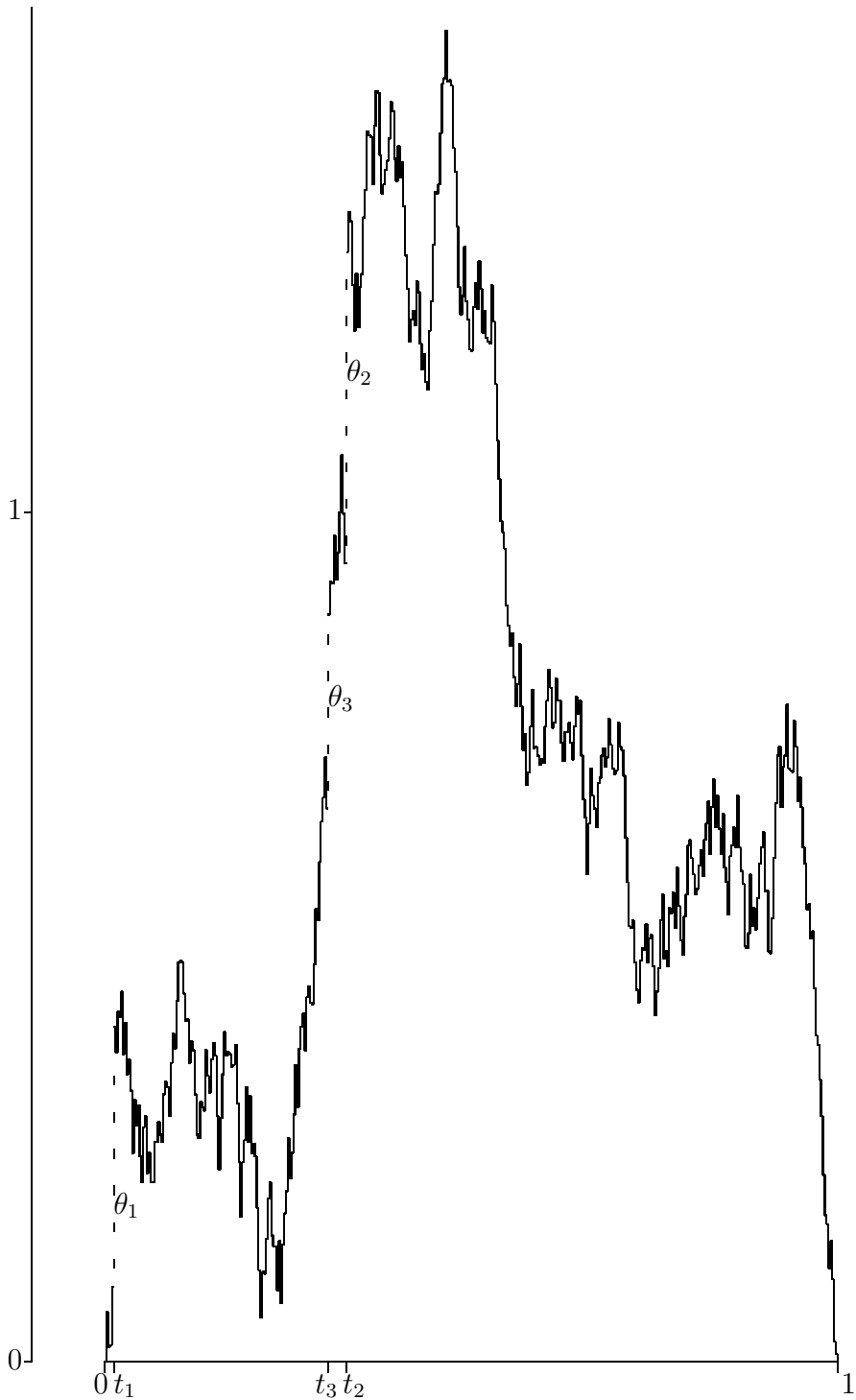}
\end{center}
\caption{A realization of $(X^{\thv}_s, 0 \leq s \leq 1)$
with $I = 3$ and
$(\theta_0,\theta_1,\theta_2,\theta_3) =
(0.862, 0.345, 0.302, 0.216)$ ($I=3$).
The jumps are marked with dashed lines; 
the jump of height $\theta_i$ occurs at time $t_i$.}
\label{fig:exc}
\end{figure}

Let $\btheta\in\bTheta$, and consider a standard Brownian
bridge $\Bbr$, and independent uniformly distributed
random variables $(U_i,i\geq1)$ in $[0,1]$, independent of $\Bbr$. Define
\begin{equation}\label{xbrtdef}
\xbrt_t=\theta_0 B^{{\rm br}}_t+\sum_{i=1}^{\infty}
\theta_i(\ind_{\{U_i\leq t\}}-t), \hskip1cm 0\leq t\leq 1.
\end{equation}
From Kallenberg \cite{kal73}, the sum on the right converges a.s.\
uniformly on $[0,1]$.
Then $\xbrt$ has exchangeable increments and infinite variation, and
by Knight \cite{knight96e} and Bertoin \cite{bertoin01eac} it
attains its overall minimum at a unique location $t_{\min}$, which is a
continuity point of $\xbrt$.
Consider the Vervaat transform $X^{\sth}$ of $\xbrt$, defined by
\begin{equation}
\lb{Xdef}
X^{\sth}_t=\xbrt_{t+t_{\min}}-\xbrt_{t_{\min}},\hskip1cm 0\leq t\leq 1,
\end{equation}
where the addition is modulo $1$. Then $X^{\sth}$ is an excursion-type process with
infinite variation, and a countable number of upward jumps with
magnitudes equal to $(\theta_i,i\geq 1)$.
See Figure \ref{fig:exc}.  Write $t_j=U_j-t_{\min}$ (mod.\ $1$)
for the location of the jump with size $\theta_j$ in $X^{\sth}$.

\begin{figure}
\begin{center}
\input{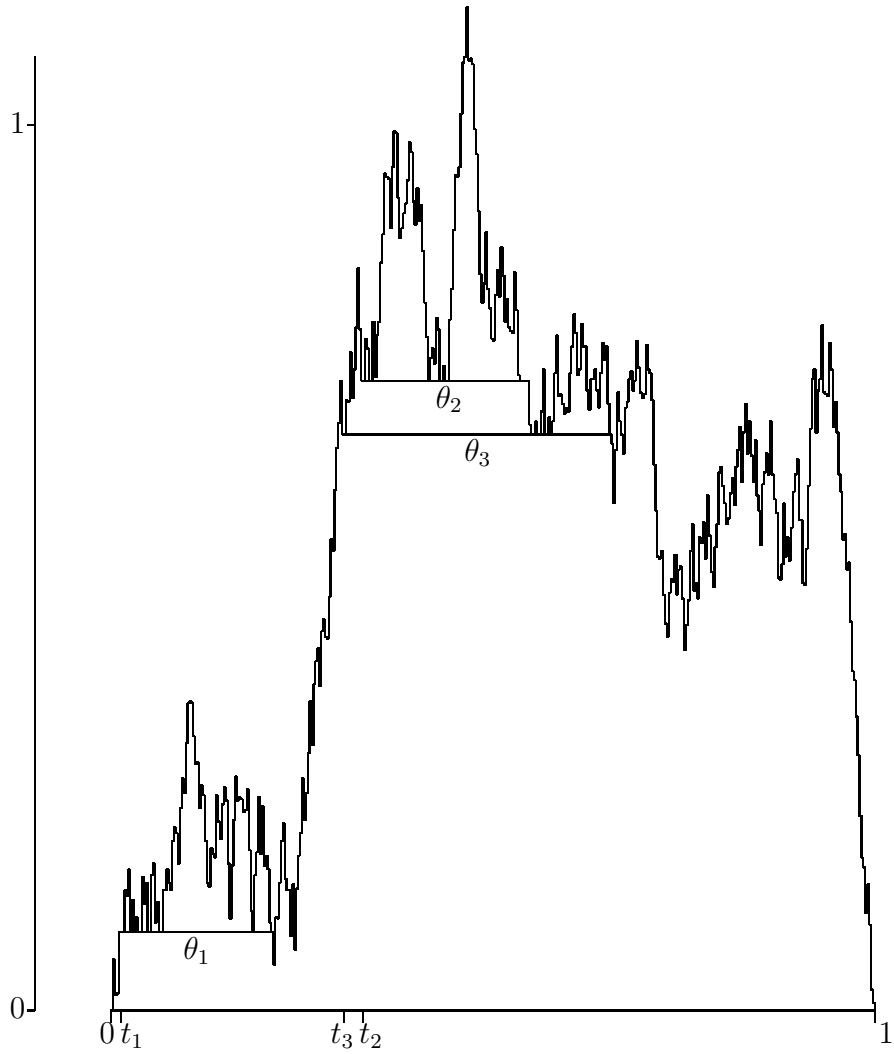}
\end{center}
\caption{The process $(Y^{\sth}_s, 0 \leq s \leq 1)$ 
constructed from the process $(X^{\sth}_s, 0 \leq s \leq 1)$
in Figure \ref{fig:exc}.
The ``reflecting'' portions of the path corresponding to jumps of $X^{\sth}$
are marked by the $\theta_i$}
\label{fig:excy}
\end{figure}

For each $j\geq 1$ such that $\theta_j>0$, write
$T_j=\inf\{s>t_j: X^{\sth}_s=X^{\sth}_{t_j-}\}$,
which exists because the process $X$ has no negative jumps. Notice
that if for some $i\neq j$ one has $t_j\in(t_i,T_i)$, then one also has
$T_j\in(t_i,T_i)$,
so the intervals $(t_i,T_i)$ are nested.
Given a sample path of $X^{\thv}$,
for $0\leq u\leq 1$ and $i\geq 1$ such that $\theta_i>0$, let
\begin{equation}
\lb{procaux}
R^{\sth}_i(u)=\left\{\begin{array}{cl}\inf_{t_{i}\leq s \leq u}X^{\sth}_{s}
-X^{\sth}_{t_i-}&\mbox{if }u\in[t_i,T_i]\\
0&\mbox{else}.
\end{array}\right.
\end{equation}
If $\theta_i=0$ then let $R^{\sth}_i$ be the null process on $[0,1]$.
We then set \begin{equation} \lb{Ydef}
Y^{\thv}=X^{\thv}-\sum_{i\geq 1}R^{\sth}_i,  \end{equation}
which is defined as the pointwise decreasing limit of
$X^{\sth}-\sum_{1\leq i\leq n}R^{\sth}_i$ as $n\to \infty$.  See Figure \ref{fig:excy}.
It is immediate that $Y^{\sth}$ is a non-negative process on $[0,1]$.
More precisely, for any $0\leq u\leq s\leq 1$ and $i$ such that
$u\geq t_i$, $R^{\sth}_i(u)$
is equal to the magnitude of the jump (if any) accomplished at time $t_i$ by
the increasing
process
$$\underleftarrow{X}^{\sth}_s(u)=\inf_{u\leq r\leq s} X^{\sth}_r,
\hskip1cm 0\leq u\leq s.$$

Since the Lebesgue measure of the range of an increasing function
$(f(s),0\leq s\leq t)$ is $f(t)-f(0)$ minus the sum of sizes of jumps accomplished by
$f$, we obtain that
\begin{equation}\label{lebesgrep}
Y^{\sth}_s=m\{\underleftarrow{X}^{\sth}_s(u):0\leq u\leq s\}\hskip1cm
0\leq s\leq 1,
\end{equation}
where $m$ is Lebesgue measure. This easily implies that $Y^{\sth}$ is a
continuous (possibly null) process, and since the largest
jump of $X^{\sth}-\sum_{1\leq i\leq n}R^{\sth}_i$ is
$\theta_{n+1}$, which tends to $0$ as $n\to \infty$, a variation of
Dini's theorem implies that
(\ref{Ydef}) holds in the sense of uniform convergence.

%



The process $Y^{\thv}$ is an
excursion-type process on $(0,1)$. Moreover, since by classical properties 
of Brownian bridges the local infima of $\xbrt$ are all distinct, the
only local
infima that $Y^{\thv}$ attains an infinite number of times are in the
intervals
$[t_i,T_i]$. Let us record some other sample path
properties of $Y^{\thv}$.

\begin{lmm}\label{XandY}
Suppose $\btheta\in
\bTheta$ has length $I\in\mathbb{N}\cup\{\infty\}$ and $\theta_0>0$.
Almost surely, the values $(X^{\sth}_{t_i-},i\geq 1)$ taken by $X^{\sth}$ at
its jump times are not attained at local minima of $X^{\sth}$. Also,
the times $t_i$ are a.s.\
not right-minima of $X^{\sth}$ in the sense that
there does not exist $\eps>0$ such that $X^{\sth}_s\geq X^{\sth}_{t_i}$
for $s\in[t_i,t_i+\eps]$.

%
\end{lmm}

\proof
Let $\xbrt_i(s)=\xbrt_s-\theta_i(\ind_{\{U_i\leq s\}}-s)$, which is independent
of $U_i$.
The shifted process $\xbrt_i(\cdot+t)-\xbrt_i(t)$ (with addition
modulo $1$) has same law as $\xbrt_i$ for every $t$, so the fact that
$1$ is not the time of a local extremum for $\xbrt_i$ and that
$|X^{{\rm br},\sth}_i(1-t)|/t\to \infty$ as $t\to 0$ (e.g.\ by
\cite[Theorem 2.2 (i)]{kal6715} and time-reversal) implies by
adding back $\theta_i(\ind_{\{U_i\leq \cdot\}}-\cdot)$ to $\xbrt_i$ that
$U_i$ is a.s.\ not a local minimum of $\xbrt$.
The statement about right-minima is obtained similarly, using the
behavior of $X^{{\rm br},\sth}_i$ at $0$ rather than $1$.

Next, since $\xbrt$ is the sum of a Brownian bridge
$\Bbr$ and an independent process,
the increments of $\xbrt$ have continuous densities, as does the Brownian
bridge (except of course the increment $\xbrt(1)-\xbrt(0)=0$ a.s.). The
probability that the minimum of $\xbrt$ in any interval $[a,b]$ with
distinct
rational bounds not containing $U_i$ equals $\xbrt_{U_i-}$ is therefore
$0$.
This finishes the proof.
\cq

The following lemma will turn out to be useful at the end of the
proof of Theorem \ref{T1}.
\begin{lmm}\label{unifconv}
Let $\btheta\in \bTheta$ satisfy $\sum_i\theta_i<\infty$, and write
$\btheta^{n}=(\theta_0,\theta_1,\ldots,\theta_n)$. Define $X^{\sth^n}$ as 
above, but where the sum defining $X^{{\rm br},\sth^n}$ is truncated at $n$.
Last, define $Y^{\sth^n}$ as in (\ref{lebesgrep}) with $X^{\sth^n}$ instead
of $X^{\sth}$. Then $Y^{\sth^n}$ converges a.s.\ uniformly to
$Y^{\sth}$ as $n\to \infty$.
\end{lmm}

\proof
We want to estimate the uniform norm $\|Y^{\sth^n}-Y^{\sth}\|$, which by
definition is $\|X^{\sth^n}-X^{\sth}-\sum_{i\geq 1}(R^{\sth^n}_i
-R^{\sth}_i)\|$ with obvious notations. The first problem is that
$X^{{\rm br},\sth^n}$ may not attain its overall infimum at the same time
as $\xbrt$, so that jump times for $X^{\sth^n}$ and $X^{\sth}$ may
not coincide anymore. So, rather than using
$X^{\sth^n}$ we consider $X'_n(s)=X^{{\rm br},\sth^n}(s+t_{\min})
-X^{{\rm br},\sth^n}(t_{\min})$ (with addition modulo $1$) where $t_{\min}$
is the time at which $\xbrt$ attains its infimum. Then $X'_n\to X^{\sth}$
uniformly. Define $R'_{n,i}$ as in (\ref{procaux}) but for the process
$X'_n$ and write $Y'_n=X'_n-\sum_{1\leq i\leq n}R'_{n,i}$. Notice that $Y'_n$
is just a slight space-time shift of $Y^{\sth^n}$, so by continuity of
$Y^{\sth^n}$ and $Y^{\sth}$ it suffices to show that $Y'_n\to Y^{\sth}$
uniformly. It is thus enough to show that
$\|\sum_{1\leq i\leq n}(R'_{n,i}-R^{\sth}_i)\|\to 0$ as $n\to \infty$.
It is easy that for each $i\geq 1$, one has uniform convergence
of $R'_{n,i}$ to $R^{\sth}_i$. Therefore, it suffices to show that
$$\lim_{k\to\infty}\limsup_{n\to\infty}\left\|
\sum_{k\leq i\leq n}R'_{n,i}\right\|
=0,$$
which is trivial because $\|R'_{n,i}\|\leq\theta_i$,
and $\sum_i \theta_i<\infty$ by hypothesis. 

\rem
Again, one guesses that the same result holds in the general $\theta_0>0$
case, so that the proof of Theorem \ref{T1} should extend to this case. 
However, the fact that $\sum_i\theta_i$ might be infinite does not
{\em a priori} prevent vanishing terms of the sum $\sum_{1\leq i\leq
n}R'_{n,i}$ to accumulate, so the proof 
might become quite technical.

\section{Constructions of {\bf p}-trees and associated excursion processes}
\label{con-ptrees}
Write $\bT_n$ for the set of rooted trees $\bt$ on vertex-set $[n]$, where
$\bt$ is directed towards its root.
Fix a probability distribution
${\bf p} = (p_1,\ldots,p_n)$.
Recall
that associated with ${\bf p}$ is a
certain distribution on $\bT_n$,
the {\em {\bf p}-tree}
\begin{equation}
P(\TT = \bt) = \prod_v p_v^{d_v},
\ d_v \mbox{ in-degree of $v$ in $\bt$} .
\lb{btree}
\end{equation}
See \cite{jp01hur} for systematic discussion of the $\bp$-tree model.
We shall define two maps
$\psi_{{\bf p}}: [0,1)^n \to \bT_n$
such that, if
$(X_1,\ldots,X_n)$ are independent $U(0,1)$ then each
$\psi_{{\bf p}}(X_1,\ldots,X_n)$
has the
distribution (\ref{btree}).
The two definitions are quite similar, but the essential difference
is that
$\psbr$ uses a {\em breadth-first} construction whereas
$\psde$ uses a {\em depth-first} construction.

\begin{figure}
\begin{center}
\setlength{\unitlength}{0.05in}
\begin{picture}(100,41)(25,36)
\put(48,37){{\footnotesize root}}
\put(50,40){a}
\put(30,50){b}
\put(50,50){c}
\put(70,50){d}
\put(20,60){e}
\put(40,60){f}
\put(70,60){g}
\put(40,70){h}
\put(49,42){\line(-2,1){16}}
\put(53,42){\line(2,1){16}}
\put(50.5,43){\line(0,1){6}}
\put(29,52){\line(-1,1){7}}
\put(32,52){\line(1,1){7}}
\put(40.5,63){\line(0,1){6}}
\put(70.5,53){\line(0,1){6}}
\put(45,70){breadth-first}

\put(108,37){{\footnotesize root}}
\put(110,40){a}
\put(90,50){b}
\put(110,50){f}
\put(130,50){h}
\put(80,60){c}
\put(100,60){d}
\put(100,70){e}
\put(110,60){g}
\put(109,42){\line(-2,1){16}}
\put(113,42){\line(2,1){16}}
\put(110.5,43){\line(0,1){6}}
\put(110.5,53){\line(0,1){6}}
\put(89,52){\line(-1,1){7}}
\put(92,52){\line(1,1){7}}
\put(100.5,63){\line(0,1){6}}
\put(105,70){depth-first}
\end{picture}
\end{center}
\caption{A planar tree, with the two orderings of vertices as
$a,b,c,d,e,f,g,h$}
\label{fig:treebfdf}
\end{figure}
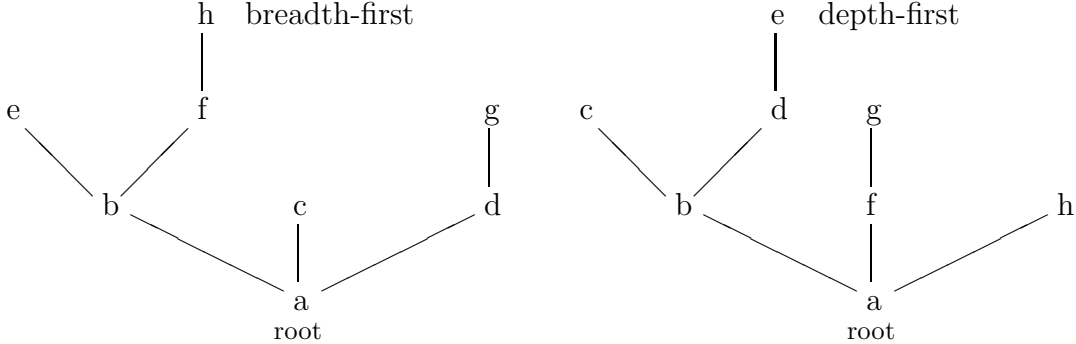

\subsection{The breadth-first construction}
The construction is illustrated in Figure \ref{fig:breadth}.
Fix distinct $(x_1,\ldots,x_n) \in [0,1)^n$.
Picture this as a configuration of particles on the circle
of unit circumference, where particle $i$ is at position $x_i$
and has a ``weight'' $p_i$ associated with it.
Define
\begin{equation}
 F^{\bp}(u) = - u +  \sum_i p_i \ind_{\{x_i \leq u\}} , \quad 0 \leq u \leq 1 . \lb{Fdef}
\end{equation}

\begin{figure}
\begin{center}
\setlength{\unitlength}{0.05in}
\begin{picture}(150,80)(-15,-5)
\put(48,37){{\footnotesize root}}
\put(50,40){4}
\put(30,50){8}
\put(50,50){2}
\put(70,50){3}
\put(20,60){7}
\put(40,60){1}
\put(70,60){5}
\put(40,70){6}
\put(49,42){\line(-2,1){16}}
\put(53,42){\line(2,1){16}}
\put(50.5,43){\line(0,1){6}}
\put(29,52){\line(-1,1){7}}
\put(32,52){\line(1,1){7}}
\put(40.5,63){\line(0,1){6}}
\put(70.5,53){\line(0,1){6}}

\put(0,0){\line(1,0){100}}
\put(0,1.1){\line(0,-1){1.1}}
\put(100,1.1){\line(0,-1){1.1}}
\put(-1,2){$0$}
\put(99,2){$1$}
\put(4,2){$x_5$}
\put(5.6,1.1){\line(0,-1){1.1}}
\put(10,-1.1){\line(0,1){1.1}}
\put(8,-5){$y(5)$}
\put(-4,-8){\vector(1,0){14}}
\put(3,-10){$p_3$}
\put(20,-5){$y(6)$}
\put(22,-1.1){\line(0,1){1.1}}
\put(12,-8){\vector(1,0){10}}
\put(15,-10){$p_7$}
\put(21,2){$x_6$}
\put(22.6,1.1){\line(0,-1){1.1}}
\put(30,-1.1){\line(0,1){1.1}}
\put(28,-5){$y(7)$}
\put(25,-10){$p_1$}
\put(24,-8){\vector(1,0){6}}
\put(39,-1.1){\line(0,1){1.1}}
\put(37,-5){$y(8)$}
\put(32,-8){\vector(1,0){7}}
\put(33,-10){$p_5$}
\put(41,-8){\vector(1,0){8}}
\put(43,-10){$p_6$}
\put(49,-1.1){\line(0,1){2.2}}
\put(47.4,2){$x_4$}
\put(47,-5){$y(1)$}
\put(54,2){$x_8$}
\put(55.6,1.1){\line(0,-1){1.1}}
\put(62,2){$x_2$}
\put(63.6,1.1){\line(0,-1){1.1}}
\put(66,2){$x_3$}
\put(67.6,1.1){\line(0,-1){1.1}}
\put(70,-1.1){\line(0,1){1.1}}
\put(68,-5){$y(2)$}
\put(51,-8){\vector(1,0){19}}
\put(58,-10){$p_4$}
\put(77,2){$x_7$}
\put(78.6,1.1){\line(0,-1){1.1}}
\put(83,2){$x_1$}
\put(84.6,1.1){\line(0,-1){1.1}}
\put(87,-1.1){\line(0,1){1.1}}
\put(85,-5){$y(3)$}
\put(72,-8){\vector(1,0){15}}
\put(79,-10){$p_8$}
\put(94,-1.1){\line(0,1){1.1}}
\put(92,-5){$y(4)$}
\put(89,-8){\vector(1,0){5}}
\put(90,-10){$p_2$}
\put(96,-8){\vector(1,0){8}}
\put(0,25){\line(1,0){100}}
\put(0,10){\line(0,1){30}}
\put(0,15){\line(-1,0){1}}
\put(0,35){\line(-1,0){1}}
\put(-5.5,34){$0.2$}
\put(-8,14){$-0.2$}
\put(0,25){\line(2,-1){5.6}}
\put(5.6,26.7){\line(2,-1){17}}
\put(22.6,23.2){\line(2,-1){26.4}}
\put(49,20.5){\line(2,-1){6.6}}
\put(55.6,25.7){\line(2,-1){8}}
\put(63.6,25.2){\line(2,-1){4}}
\put(67.6,31.2){\line(2,-1){11}}
\put(78.6,31.7){\line(2,-1){6}}
\put(84.6,32.7){\line(2,-1){15.4}}
\put(10,30){$F^{\bp}(u)$}
\end{picture}
\end{center}
\caption{The construction of the tree $\psbr (x_1,\ldots,x_8)$}
\label{fig:breadth}
\end{figure}

There exists some particle $v$ such that $F^{\bp}(x_v-)=
\inf_u F^{\bp}(u)$: assume the particle is unique.
Let $v=v_1,v_2,\ldots,v_n$ be the ordering of particles
according to the natural ordering of positions
$x_{v_1} < x_{v_2} < \ldots$
around the circumference of the circle.
(In Figure \ref{fig:breadth} we have $v_1 = 4$ and the ordering is $4, 8, 2, 3, 7, 1, 5, 6$).
Write $y(1) = x_{v_1}$ and for $2 \leq j \leq n$ let
$y(j+1) = y(j) + p_{v_j} \bmod 1$.
So $y(n+1) = y(1)$ and the successive intervals $[y(j),y(j+1)], 1 \leq j \leq n$
are adjacent and cover the circle.
We assert
\begin{equation}
x_{v_j} \in [y(1),y(j)), \ \ 2 \leq j \leq n .
\lb{claim}
\end{equation}
To argue by contradiction,
suppose this fails first for $j$.
Then $[y(1),y(j))$, interpreted $\bmod\ 1$, contains particles
$v_1,\ldots,v_{j-1}$ only.
Since $y(j)-y(1) = p_{v_1} + \ldots + p_{v_{j-1}}$ this implies
$F^{\bp}(y(j)-) = F^{\bp}(y(1)-)$,
contradicting uniqueness of the minimum.

We specify the tree $\psbr (x_1,\ldots,x_n)$ by:
\begin{quote}
        $v_1$ is the root\\
the children of $v_j$ are the particles $v$ with
$x_v \in (y(j),y(j+1))$.
\end{quote}
By (\ref{claim}), any child $v_k$ of $v_j$ has $k>j$,
so the graph cannot contain a cycle.
If it were a forest and not a single tree, then the component containing
the root $v_1$ would consist of vertices $v_1,\ldots, v_j$ for some $j<n$.
Then the interval $[y(1),y(j+1)]$ would contain only the
particles $v_1,\ldots,v_j$,
contradicting (\ref{claim}) for $j+1$.

Thus the construction does indeed give a tree.
From the viewpoint of this construction it would be natural
to regard the tree as {\em planar}
(or {\em ordered}:
the $d_v$ children of $v$ are distinguished as first, second, etc)
but we disregard order and view trees in $\bT_n$ as unordered.

Now consider the case where
$(x_1,\ldots,x_n)=(X_1,\ldots,X_n)$ are independent $U(0,1)$.
Fix an unordered tree $\bt$ and write $v_1$ for its root.
Fix an arbitrary $x_{v_1} \in (0,1)$ and condition on $X_{v_1} = x_{v_1}$.
Consider the chance that the construction yields the particular tree $\bt$.
For this to happen,
the particles corresponding to the $d_{v_1}$ children of $v_1$
must fall into the interval $[x_{v_1},x_{v_1} + p_{v_1}]$,
which has chance $p_{v_1}^{d_{v_1}}$.
Inductively, for each vertex $v$ an interval of length $p_v$ is specified
and it is required that $d_v$ specified particles fall into that interval,
which has chance $p_v^{d_v}$.
So the conditional probability of constructing $\bt$ is indeed the probability in (\ref{btree}),
and hence so is the unconditional probability.

\rem Note that in the argument above we do not start by conditioning
on $F^{\bp}$ having its minimum at $x_{v_1}$, which would affect the
distribution of the $(X_i)$.

We now derive an interpretation (\ref{F1},\ref{F2})
of the function $F^{\bp}$ at (\ref{Fdef}),
which will be used in the asymptotic setting later. From now on we
also suppose that for $j\geq 2$, $y(j)$ is not a jump time for
$F^{\bp}$ to avoid needing the distinction between
$F^{\bp}(y(j))$ and $F^{\bp}(y(j)-)$;
this is obviously true a.s.\ when the jump times
are independent uniform, which will be the relevant case.

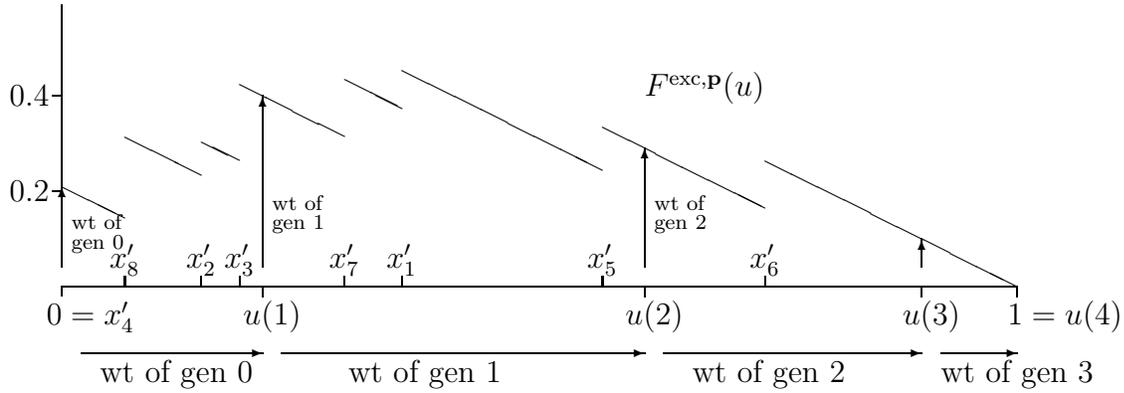
\begin{figure}
\begin{center}
\setlength{\unitlength}{0.05in}
\begin{picture}(100,50)(55,-14)
\put(49,0){\line(1,0){100}}
\put(149,-1.1){\line(0,1){1.1}}
\put(148,-4){$1 = u(4)$}
\put(104,2){$x^\prime_5$}
\put(105.6,1.1){\line(0,-1){1.1}}
\put(110,-1.1){\line(0,1){1.1}}
\put(108,-4){$u(2)$}
\put(110,2){\vector(0,1){12.5}}
\put(111,8){{\scriptsize wt of}}
\put(111,6){{\scriptsize gen 2}}
\put(112,-7){\vector(1,0){27}}
\put(115,-10){wt of gen 2}
\put(121,2){$x^\prime_6$}
\put(122.6,1.1){\line(0,-1){1.1}}
\put(139,-1.1){\line(0,1){1.1}}
\put(139,2){\vector(0,1){3}}
\put(137,-4){$u(3)$}
\put(141,-7){\vector(1,0){8}}
\put(141,-10){wt of gen 3}
\put(47.4,-4){$0 =x^\prime_4$}
\put(54,2){$x^\prime_8$}
\put(55.6,1.1){\line(0,-1){1.1}}
\put(62,2){$x^\prime_2$}
\put(63.6,1.1){\line(0,-1){1.1}}
\put(66,2){$x^\prime_3$}
\put(67.6,1.1){\line(0,-1){1.1}}
\put(70,-1.1){\line(0,1){1.1}}
\put(70,2){\vector(0,1){18}}
\put(71,8){{\scriptsize wt of}}
\put(71,6){{\scriptsize gen 1}}
\put(68,-4){$u(1)$}
\put(51,-7){\vector(1,0){19}}
\put(53,-10){wt of gen 0}
\put(77,2){$x^\prime_7$}
\put(78.6,1.1){\line(0,-1){1.1}}
\put(83,2){$x^\prime_1$}
\put(84.6,1.1){\line(0,-1){1.1}}
\put(72,-7){\vector(1,0){38}}
\put(79,-10){wt of gen 1}
\put(49,0){\line(1,0){100}}
\put(49,-1){\line(0,1){1}}
\put(49,10.5){\line(0,1){19}}
\put(49,2){\vector(0,1){8.5}}
\put(50,6){{\scriptsize wt of}}
\put(50,4){{\scriptsize gen 0}}
\put(49,10){\line(-1,0){1}}
\put(49,20){\line(-1,0){1}}
\put(43.5,9){$0.2$}
\put(43.5,19){$0.4$}
\put(105.6,16.7){\line(2,-1){17}}
\put(122.6,13.2){\line(2,-1){26.4}}
\put(49,10.5){\line(2,-1){6.6}}
\put(55.6,15.7){\line(2,-1){8}}
\put(63.6,15.2){\line(2,-1){4}}
\put(67.6,21.2){\line(2,-1){11}}
\put(78.6,21.7){\line(2,-1){6}}
\put(84.6,22.7){\line(2,-1){21}}
\put(110,20){$\Fexc(u)$}
\end{picture}
\end{center}
\caption{$\Fexc(\cdot)$ codes the weights of successive generations
({\em wt}  of {\em gen}) of the $\bp$-tree in Figure \ref{fig:breadth}}
\label{fig:fexcp}
\end{figure}

For $2 \leq j \leq n$, vertex $v_j$ has some parent
$v_{z(j)}$, where $1 \leq z(j) < j$.
By induction on $j$,
\[ F^{\bp}(y(j)) - F^{\bp}(y(1)-)
= \sum_{i: i>j, z(i) \leq j} p_{v_i}
. \]
In words, regarding $\bt$ as ordered, the sum is over vertices
$i$ which are in the same generation as $j$ but {\em later} than $j$;
and over vertices $i$ in the next generation whose parents are
{\em before} $j$ or are $j$ itself.
For $h \geq 1$,
write $t(h)$ for the number of vertices at height $\leq h-1$.
The identity above implies
\[ F^{\bp}(y(t(h)+1)) - F^{\bp}(y(1)-)
= \sum_{v: \HT(v) =h} p_v . \]
Also by construction
\[ y(t(h)+1) - y(1) \bmod 1
= \sum_{v: \HT(v) \leq h-1} p_v . \]
We can rephrase the last two inequalities in terms of the ``excursion'' function
\begin{equation}
\lb{Fexc}
 \Fexc(u) := F^{\bp}(y(1) + u \bmod 1)-F^{\bp}(y(1)-), \ 0 \leq u \leq 1
\end{equation}
and of
$u(h) := y(t(h)+1) - y(1) \bmod 1$.
Then
\begin{eqnarray}
u(h) &=&
 \sum_{v: \HT(v) \leq h-1} p_v  \lb{F1} \\
\Fexc(u(h)) &=&
 \sum_{v: \HT(v) =h} p_v . \lb{F2}
\end{eqnarray}
So the weights of successive generations are coded within $\Fexc(\cdot)$,
as illustrated in
Figure \ref{fig:fexcp}. 
Note
that to draw Figure \ref{fig:fexcp} we replace $x_i$ by
\[ x^\prime_i := x_i - y(1) \bmod 1 . \]

\rem
There is a queuing system interpretation to the breadth-first
construction, which was pointed out to us by a referee. 
In this interpretation, the customer labelled $i$ arrives at time $x'_i$ and
requires a total service time $p_i$.  If customers are served according to the
FIFO rule (first-in first-out) then $\fexp(u)$ is the remaining amount
of time needed to serve the customers in line at time $u$.

\subsection{The depth-first construction}\label{subdfc}

The construction is illustrated in Figure \ref{fig:depth},
using the same $(x_i)$ and $(p_i)$ as before, and hence
the same $F^{\bp}(u)$.
In the previous construction we ``examined'' particles in the order
$v_1,v_2,\ldots,v_n$;
we defined $y(1) = v_1$ and inductively\\
$\bullet$ $y(j+1) = y(j) + p_{v_j} \bmod 1$\\
$\bullet$
the children of $v_j$ are the particles $v$ with
$x_v \in (y(j),y(j+1))$.
\\
In the present construction we shall examine particles in a different order
$w_1,w_2,\ldots,w_n$
and use different $y^\prime(j)$
to specify the intervals which determine the offspring of a parent.
Start as before with $w_1 = v_1$ and $y^\prime(1) = x_{w_1}$.
Inductively set\\
$\bullet$ $y^\prime(j+1) = y^\prime(j) + p_{w_j} \bmod 1$\\
$\bullet$
the children of $w_j$ are the particles $v$ with
$x_v \in (y^\prime(j),y^\prime(j+1))$.\\
$\bullet$ $w_{j+1}$ is

the first child of $w_j$, if any; else

the next unexamined child of parent($w_j$), if any; else

the next unexamined child of parent(parent($w_j$)), if any; else

and so on.\\

\noindent
Here ``unexamined'' means ``not one of $w_1,\ldots,w_j$'' and ``next'' uses the
natural order of children of the same parent.

Figure \ref{fig:depth} and the following paragraph talk through the 
construction in a particular example,
using the same $(x_i)$ and $(p_i)$ as in Figure \ref{fig:breadth}. 
Checking that $\psde (X_1,\ldots,X_n)$ has distribution (\ref{btree}),
i.e. is a random $\bp$-tree, uses exactly the same argument as before.

\begin{figure}
\begin{center}
\setlength{\unitlength}{0.05in}
\begin{picture}(100,94)(5,-14)
\put(48,37){{\footnotesize root}}
\put(50,40){4}
\put(30,50){8}
\put(50,50){2}
\put(70,50){3}
\put(20,60){7}
\put(40,60){1}
\put(40,70){5}
\put(50,60){6}
\put(49,42){\line(-2,1){16}}
\put(53,42){\line(2,1){16}}
\put(50.5,43){\line(0,1){6}}
\put(50.5,53){\line(0,1){6}}
\put(29,52){\line(-1,1){7}}
\put(32,52){\line(1,1){7}}
\put(40.5,63){\line(0,1){6}}

\put(0,0){\line(1,0){100}}
\put(0,1.1){\line(0,-1){1.1}}
\put(100,1.1){\line(0,-1){1.1}}
\put(-1,2){$0$}
\put(99,2){$1$}
\put(4,2){$x_5$}
\put(5.6,1.1){\line(0,-1){1.1}}
\put(7,-1.1){\line(0,1){1.1}}
\put(5,-5){$y^\prime(5)$}
\put(-4,-8){\vector(1,0){11}}
\put(1,-10){$p_1$}
\put(14,-5){$y^\prime(6)$}
\put(16,-1.1){\line(0,1){1.1}}
\put(9,-8){\vector(1,0){7}}
\put(11,-10){$p_5$}
\put(21,2){$x_6$}
\put(22.6,1.1){\line(0,-1){1.1}}
\put(23,-1.1){\line(0,1){1.1}}
\put(21,-5){$y^\prime(7)$}
\put(19,-10){$p_2$}
\put(18,-8){\vector(1,0){5}}
\put(33,-1.1){\line(0,1){1.1}}
\put(31,-5){$y^\prime(8)$}
\put(25,-8){\vector(1,0){8}}
\put(28,-10){$p_6$}
\put(35,-8){\vector(1,0){14}}
\put(41,-10){$p_3$}
\put(49,-1.1){\line(0,1){2.2}}
\put(47.4,2){$x_4$}
\put(47,-5){$y^\prime(1)$}
\put(54,2){$x_8$}
\put(55.6,1.1){\line(0,-1){1.1}}
\put(62,2){$x_2$}
\put(63.6,1.1){\line(0,-1){1.1}}
\put(66,2){$x_3$}
\put(67.6,1.1){\line(0,-1){1.1}}
\put(70,-1.1){\line(0,1){1.1}}
\put(68,-5){$y^\prime(2)$}
\put(51,-8){\vector(1,0){19}}
\put(58,-10){$p_4$}
\put(77,2){$x_7$}
\put(78.6,1.1){\line(0,-1){1.1}}
\put(83,2){$x_1$}
\put(84.6,1.1){\line(0,-1){1.1}}
\put(87,-1.1){\line(0,1){1.1}}
\put(85,-5){$y^\prime(3)$}
\put(72,-8){\vector(1,0){15}}
\put(79,-10){$p_8$}
\put(99,-1.1){\line(0,1){1.1}}
\put(97,-5){$y^\prime(4)$}
\put(89,-8){\vector(1,0){10}}
\put(93,-10){$p_7$}
\put(101,-8){\vector(1,0){3}}
\put(0,25){\line(1,0){100}}
\put(0,10){\line(0,1){30}}
\put(0,15){\line(-1,0){1}}
\put(0,35){\line(-1,0){1}}
\put(-5.5,34){$0.2$}
\put(-8,14){$-0.2$}
\put(0,25){\line(2,-1){5.6}}
\put(5.6,26.7){\line(2,-1){17}}
\put(22.6,23.2){\line(2,-1){26.4}}
\put(49,20.5){\line(2,-1){6.6}}
\put(55.6,25.7){\line(2,-1){8}}
\put(63.6,25.2){\line(2,-1){4}}
\put(67.6,31.2){\line(2,-1){11}}
\put(78.6,31.7){\line(2,-1){6}}
\put(84.6,32.7){\line(2,-1){15.4}}
\put(10,30){$F^{\bp}(u)$}
\end{picture}
\end{center}
\caption{The construction of the tree $\psde (x_1,\ldots,x_n)$}
\label{fig:depth}
\end{figure}

As in Figure \ref{fig:breadth}, the root of the tree is vertex $4$ ($w_1 = 4$),
and we set
$y^\prime(1) = x_{4}$.
As before, $y^\prime(2) = y^\prime(1) + p_{4}$, and the children of the root are
the vertices $\{8,2,3\}$
for which $x_v \in (y^\prime(1),y^\prime(2))$.
As before, we next examine the first child $w_2 = 8$ of the root,
set $y^\prime(3) = y^\prime(2) + p_{8}$, and let the children of $8$ be
the vertices $\{7,1\}$
for which $x_v \in (y^\prime(2),y^\prime(3))$.
At this stage the constructions differ.
We next examine vertex $7$, being the first child of vertex $8$,
by setting $y^\prime(4) = y^\prime(3)+p_7$; the children of vertex $8$
are the vertices $v$ with $x_v \in (y^\prime(3),y^\prime(4))$, and it turns
out there are no such vertices.
We continue examining vertices in the depth-first order
$4,8,7,1,5,2,6,3$.

As with the breadth-first construction, the point of the
depth-first construction is that the excursion function
$\Fexc(\cdot)$ tells us something about the distribution
of the tree.
For each vertex $v$ of
$\psde (x_1,\ldots,x_n)$
there is a path
${\rm root} = y_0,y_1,\ldots,y_j = v$
from the root to $v$.
For each $0 \leq i < j$ the vertex $y_{i+1}$ is a child of vertex $y_i$;
let $y_{i,1},y_{i,2},\ldots$ be the later children of $y_i$,
and let $y_{j,1},y_{j,2},\ldots$ be all children of $v$.
Write
${\cal N}(v) = \cup_{0 \leq i \leq j} \{y_{i,1},y_{i,2},\ldots\}$.

\begin{figure}
\begin{center}
\setlength{\unitlength}{0.05in}
\begin{picture}(100,55)(55,-13)
\put(49,0){\line(1,0){100}}
\put(49,1.1){\line(0,-1){1.1}}
\put(149,-1.1){\line(0,1){1.1}}
\put(149,-4){$1$}
\put(104,2){$x^\prime_5$}
\put(105.6,1.1){\line(0,-1){1.1}}
\put(110,-1.1){\line(0,1){1.1}}
\put(121.3,2){$x^\prime_6$}
\put(122.6,1.1){\line(0,-1){1.1}}
\put(116,-1.1){\line(0,1){1.1}}
\put(114,-4){$y^*(6)$}
\put(123,-1.1){\line(0,1){1.1}}
\put(121,-4){$y^*(7)$}
\put(116,-7){\vector(1,0){7}}
\put(114,-10){examine $2$}
\put(123,2){\vector(0,1){11}}
\put(119,5.3){$p_6+p_3$}
\put(47.4,-4){$0 =x^\prime_4$}
\put(54,2){$x^\prime_8$}
\put(55.6,1.1){\line(0,-1){1.1}}
\put(62,2){$x^\prime_2$}
\put(63.6,1.1){\line(0,-1){1.1}}
\put(66,2){$x^\prime_3$}
\put(67.6,1.1){\line(0,-1){1.1}}
\put(70,-1.1){\line(0,1){1.1}}
\put(68,-4){$y^*(1)$}
\put(77,2){$x^\prime_7$}
\put(78.6,1.1){\line(0,-1){1.1}}
\put(83,2){$x^\prime_1$}
\put(84.6,1.1){\line(0,-1){1.1}}
\put(87,-1.1){\line(0,1){1.1}}
\put(85,-4){$y^*(2)$}
\put(87,2){\vector(0,1){19.5}}
\put(77.3,10.7){$p_7+p_1+p_2+p_3$}
\put(70,-7){\vector(1,0){17}}
\put(71,-10){examine $8$}
\put(49,0){\line(1,0){100}}
\put(49,-1){\line(0,1){30}}
\put(49,10){\line(-1,0){1}}
\put(49,20){\line(-1,0){1}}
\put(43.5,9){$0.2$}
\put(43.5,19){$0.4$}
\put(105.6,16.7){\line(2,-1){17}}
\put(122.6,13.2){\line(2,-1){26.4}}
\put(49,10.5){\line(2,-1){6.6}}
\put(55.6,15.7){\line(2,-1){8}}
\put(63.6,15.2){\line(2,-1){4}}
\put(67.6,21.2){\line(2,-1){11}}
\put(78.6,21.7){\line(2,-1){6}}
\put(84.6,22.7){\line(2,-1){21}}
\put(110,20){$\Fexc(u)$}
\end{picture}
\end{center}
\caption{Relation \re{deprop} in the depth-first construction}
\label{fig:deprop}
\end{figure}
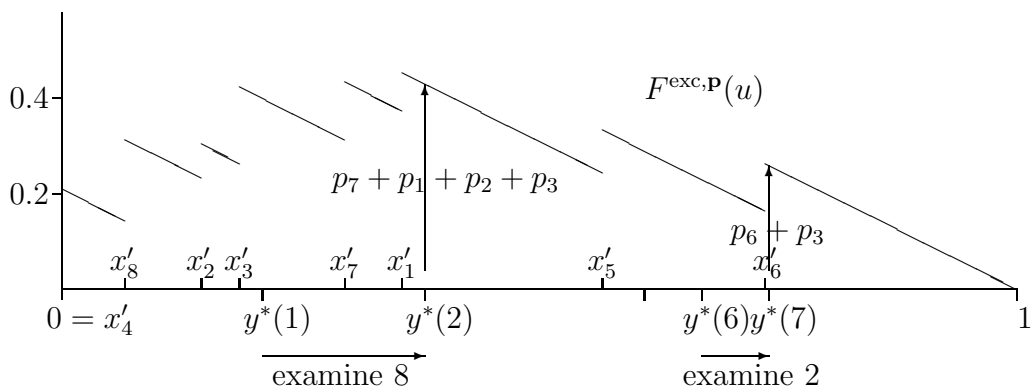

In the $u$-scale of $\Fexc(u)$, we finish ``examining''
vertex $w_i$ at time
$y^*(i):=y^\prime(i) - y^\prime (1)$.
For vertex $v = w_i$ set $e(v) = y^*(i)$.
Then the relevant property of $\Fexc$ is
\begin{equation}
\Fexc(e(v)) = \sum_{w \in {\cal N}(v)} p_w, \quad
\forall v .
\lb{deprop}
\end{equation}
See Figure \ref{fig:deprop} for illustration.
As before, in Figure \ref{fig:deprop} the position of the jump of height $p_i$
is moved from $x_i$ to
$x^\prime_i := x_i - y^\prime(1) \bmod 1$.
At first sight, relation \re{deprop}
may not look useful.
But we shall see in section \ref{skocon}
that in the asymptotic regime the right
side of \re{deprop} can be related to
$\sum_{w \mbox{ ancestor of } v} p_w$
which in turn relates to the height of $v$.

\rem
We might alternatively have defined the tree
$\psde(x_1,\ldots,x_n)$ in a way that would have been less suited
for the forthcoming analysis, but which is worth mentioning. It is based on
the LIFO-queuing system construction of Galton-Watson trees in
Le Gall-Le Jan \cite{lglj97} which we sketch here.
Imagine vertex $i$ is a customer in a line which requires a
treatment time $p_i$. The customer $i$ arrives at time $x_i$ and customers are
treated according to the Last In First Out rule. After relocating the 
the time-origin is at the time when the minimum of the bridge $F^{\bp}$ is 
attained, 
the first customer in line will also be the last to get out. Then we say that
vertex $i$ is a parent of vertex $j$ if customer $j$ arrives in a time-interval
when $i$ was being treated. Notice that the tree thus defined is in general
different from $\psde(x_1,\ldots,x_n)$.

It is easy to see, using induction and the same kind of
arguments as above, that taking $x_1,\ldots,x_n$ to be independent uniform
random variables builds a $\bp$-tree (in order that $i$ has $k$ children, $k$
uniform random variables must interrupt the service of $i$ which takes total
time $p_i$, so this has probability $p_i^k$). It is also easy that the order
of customer arrivals (after relocating the time origin) corresponds to the
depth-first order on the tree. In particular, the cyclic depth-first
random order of vertices in a $\bp$-tree is the uniform cyclic order on the
$n$ vertices.

\section{Convergence of $\bp$-trees to the ICRT}\label{ptree-icrt}
Here we review known results concerning
convergence of $\bp$-trees to the ICRT,
and spotlight what new results are required to prove Theorems \ref{T1} and
\ref{T2}.

The general notion \re{def-explore} of exploration process of a continuum
random tree can be reinterpreted as follows.
Fix $J \geq 1$.
Let $(U_j, 1 \leq j \leq J)$ be independent $U(0,1)$ r.v.s and let
$U_{(1)}<U_{(2)} < \ldots < U_{(J)}$
be their order statistics.
To an excursion-type process
$(H_s, 0 \leq s \leq 1)$
associate the random $2J-1$-vector
\begin{equation}
\lb{ZUJ}
 \left( H_{U_{(1)}}, \inf_{U_{(1)}\leq s \leq U_{(2)}} H_s,
  H_{U_{(2)}}, \inf_{U_{(2)}\leq s \leq U_{(3)}} H_s,
\ldots, H_{U_{(J)}} \right) .
\end{equation}
This specifies a random tree-with-edge-lengths, with $J$ leaves,
as follows.\\
$\bullet$ The path from the root to the $i$'th leaf has length
$H_{U_{(i)}}$.\\
$\bullet$ The paths from the root to the $i$'th leaf
and from the root to the $(i+1)$'st leaf have their branchpoint
at distance $\inf_{U_{(i)} \leq s \leq U_{(i+1)}} H_s$.\\
Now label the $i$'th leaf as vertex $i^\prime$,
where $U_{(i)} = U_{i^\prime}$. 
Write the resulting tree as $\TT^H_J$. 
Call this the {\em sampling a function} construction. 

On the  other hand one  can use a  continuum random tree $\TT$
to define a random tree-with-edge-lengths $\TT_J$ as follows.\\
$\bullet$ Take a realization of $\TT$.\\
$\bullet$ From the mass measure on that realization, pick independently
$J$ points and label them as $\{1,2,\ldots,J\}$.\\
$\bullet$ Construct the spanning tree on those $J$ points and the root; 
this is the realization of $\TT_J$.

\noindent
Call this the {\em sampling a CRT} construction.

As discussed in detail in \cite{me56},
the relationship
\begin{quote}
the exploration process of $\TT$ is distributed as $(H_s, 0 \leq s \leq 1)$
\end{quote}
is equivalent to
\[
\TT_J \ed \TT^H_J, \quad \forall \ J \geq 1 ,
\]
(the background hypotheses in \cite{me56} were rather different,
assuming path-continuity for instance, but the ideas go through
to our setting.) 
In our setting, there is an explicit description of
the distribution of the spanning tree
$\TT^{\sth}_J$ derived from the ICRT $\TT^{\sth}$
(see \cite{me88}), so to prove Theorem \ref{T1}
it is enough to verify
\begin{equation}
\lb{TTZJ}
\TT^{\sth}_J \ed \TT^{2Y/\theta_0^2}_J, \quad \forall \ J \geq 1
\end{equation}
for $Y = Y^{\sth}$ defined at \re{Ydef}.
In principle one might verify \re{TTZJ} directly, but this
seems difficult even in the
case $J = 1$.
Instead we shall rely on weak convergence arguments, starting
with the known Proposition \ref{PwcICRT} below.

Consider a probability distribution
$\bp = (p_1,\ldots,p_n)$ which is {\em ranked}:
$p_1 \geq p_2 \geq \ldots \geq p_n > 0$.
In the associated $\bp$-tree \re{btree}, pick $J$ vertices independently
from distribution $\bp$, label them as
$[J]$ in order of pick, take the spanning tree on the root and these $J$ vertices,
regard each edge as having length $1$,
and then delete degree-$2$ vertices to form edges of positive integer length. 
Call the resulting random tree ${\cal S}^{\bp}_J$.
Define
$\sp := \sqrt{\sum_i p_i^2}$.
Now consider a sequence $\bp_n = (p_{ni})$ of ranked probability
distributions which satisfy
\begin{equation}
\lb{regime}
\lim_n \sigma(\bp_n) =0; \quad
\lim_n p_{ni}/\sigma(\bp_n) = \theta_i , \ 1 \leq i \leq I ;
\quad \lim_n p_{ni}/\sigma(\bp_n) = 0, \ i>I
\end{equation}
for some limit $\btheta = (\theta_0, \ldots,\theta_I) \in \Thfinite$. 
For a tree $\bt$ and a real constant $\sigma > 0$ define
$\sigma \otimes \bt$ to be the tree obtained from $\bt$ by
multiplying edge-lengths by $\sigma$.
The following result summarizes Propositions 2, 3 and 5(b) of \cite{me87}.
Recall
$\TT^{\sth}_J$ is obtained by sampling the ICRT $\TT^{\sth}$.
\begin{prp}\label{PwcICRT}
For a sequence $\bp = \bp_n$ satisfying \re{regime},
as $n \to \infty$
\[ \sp \otimes {\cal S}^\bp_J \cd \TT^{\sth}_J, \ J \geq 1 . \]
\end{prp}
The tree ${\cal S}^\bp_J$ may not be well-defined because
two of the $J$ sampled vertices may be the same; but part of
Proposition \ref{PwcICRT} is that this probability tends to zero.

Now consider the ``bridge'' process $F^\bp$ at \re{Fdef}, where from
now on the jump times $x_1,\ldots,x_n$ are uniformly distributed independent
random variables.
Standard results going back to Kallenberg \cite{kal73}
show that, under the asymptotic regime \re{regime},
\[ (\sigma^{-1}(\bp) F^{\bp}(s), \ 0 \leq s \leq 1) \cd
( \xbrt_s, \ 0 \leq s \leq 1  ) , \]
where $\xbrt$ is defined at \re{xbrtdef}.
It follows by an argument that can be found e.g.\ in \cite{bertoin01eac}
(using the continuity of the bridge process at its minimum)
that the associated excursion process
$\fexp$ at \re{Fexc}
satisfies
\begin{equation}
\lb{pFX}
(\sigma^{-1}(\bp) \fexp(s), \ 0 \leq s \leq 1)
\cd
(X^{\sth}_s, \ 0 \leq s \leq 1)
\end{equation}
for $X^\sth$ defined at \re{Xdef}.  

Recall from section \ref{sec-cXZ}
how $(Y_s^{\thv})$ is constructed as a modification of $(X_s^{\thv})$.
We next describe a parallel modification of $\fexp$
to construct a process $G^{\bp}_I$.
Given a realization of the $\bp$-tree obtained via the
depth-first construction illustrated in Figure \ref{fig:deprop},
and given
$I \geq 0$,
let $\BB_{i} \subseteq [n]$
be the set of vertices which are the child of some vertex $i$ in
from $\{1,\ldots,I\}$.
In the setting of the depth-first construction of the
$\bp$-tree from $\fexp$,
illustrated in Figure \ref{fig:deprop}, for every vertex $v 
\in\BB_{i}$, define
\begin{equation}\lb{rhok}
\begin{array}{ccll}
\rho_v(u) &=& 0&  0 \leq u \leq x^\prime_i\\
&=& p_v&  x^\prime_i < u \leq e(v)-p_v\\
&=& e(v) - u&  e(v)-p_v \leq u \leq e(v)\\
&=& 0&  e(v) \leq u \leq 1 .
\end{array}
\end{equation}
and then let
\begin{equation}
\lb{ri}
r_i^{\bp}(u)=\sum_{v\in\BB_i}\rho_v(u)
\end{equation}
and
\begin{equation}
\lb{GIp}
G^{\bp}_I(u) = \fexp(u) - \sum_{i=1}^I r_i^{\bp}(u) .
\end{equation}
We will show in section \ref{convY} that
\re{pFX} extends to
\begin{prp}\label{pGZ}
For a sequence $\bp = \bp_n$ satisfying \re{regime}
with limit $\btheta = (\theta_0,\ldots,\theta_I)$,
as $n \to \infty$
\[
(\sigma^{-1}(\bp) G^{\bp}_I(u), \ 0 \leq u \leq 1)
\cd
(Y^{\sth}(u), \ 0 \leq u \leq 1)
\]
for $Y^\sth$ defined at \re{Ydef}, for the Skorokhod topology.
\end{prp}


We finally come to the key issue; we want to show that
$G^{\bp}_I(\cdot)$
approximates the (discrete) exploration process.
In the depth-first construction of the $\bp$-tree $\TT$ from
$\fexp$,
we examine vertex $w_i$ during $(y^*(i-1),y^*(i)]$.
Define
\begin{equation}
 H^\bp (u) :=
\mbox{ height of $w_i$ in $\TT$ ;}
\quad u \in (y^*(i-1),y^*(i)] . \lb{defHp} 
\end{equation}
Roughly, we show
that {\em realizations} of
$\frac{\theta_0^2}{2} \sp
H^\bp(\cdot)$
and of
$\sigma^{-1}(\bp) G^\bp_I(\cdot)$
are close.
Precisely, we will prove the following in section \ref{skocon}
\begin{prp}\label{Pkey}
Let $\btheta \in \Thfinite$.
There exists a sequence $\bp = \bp_n$ satisfying \re{regime}
with limit $\btheta$,
such that as $n \to \infty$,
$$\sup_{u\in[0,1]}\left| \sfrac{\theta_0^2}{2}\sp H^{\bp}(u)-
\sp^{-1}\gip(u) \right|\cp  0. $$
\end{prp}

The next result, Lemma \ref{Lrelate}, relates
the exploration process $H^\bp$ at \re{defHp}
to the spanning trees
${\cal S}^{\bp}_J$.
This idea was used in (\cite{me98}; proof of Proposition 7)
but we say it more carefully here.
Given $u^1 \in (0,1)$ define, as in \re{defHp},
\[ w^1 = w_i \mbox{ for $i$ specified by }
u^1 \in (y^*(i-1),y^*(i)] . \]
Given $0<u^1<u^2<1$, define $w^2$ similarly, and let vertex $b$
be the branchpoint of the paths from the root to vertices $w^1$
and $w^2$.
Distinguish two cases.\\
{\em Case (i)}:
$w^1=w^2$ or $w^1$ is an ancestor of $w^2$.
In this case $b = w^1$ and so trivially
$\HT(b) = \min_{u^1\leq u \leq u^2} H^{\bp}(u)$.\\
{\em Case (ii)}:
otherwise, $b$ is a strict ancestor of both $w^1$ and $w^2$.
In this case we assert
\[ \HT(b) = \min_{u^1\leq u \leq u^2} H^{\bp}(u) - 1,\]
because vertex $b$ appears, in the depth-first order, strictly before
vertex $w^1$.  Then consider the set of vertices
between $w^1$ and $w^2$ (inclusive) in the depth first order.
This set contains the child $w^*$ of $b$ which is an ancestor of
$w^2$ or is $w^2$ itself, and $\HT(w^*) = \HT(b) +1$.
But the set cannot contain any vertex of lesser height.

Now the length of the interval $(y^*(i-1),y^*(i)]$ equals $p_{w_i}$
by construction.
So if $U^1$ has uniform distribution on $(0,1)$ then the
corresponding vertex $W^1$ at \re{defHp} has distribution $\bp$.
Combining with the discussion above regarding branchpoint heights gives
\begin{lmm}\label{Lrelate}
Fix $\bp$, make the depth-first construction of a $\bp$-tree
and define $H^{\bp}$ by \re{defHp}.
Fix $J$.
Take $U_1,\ldots,U_J$ independent uniform $(0,1)$
and use them and $H^{\bp}$ to define a tree-with-edge-lengths
$\TT^{\bp}_J$ via the ``sampling a function''
construction below \re{ZUJ}.
Then this tree agrees, up to perhaps changing heights of branchpoints
by $1$, with a tree distributed as the tree
${\cal S}^{\bp}_J$ defined above \re{regime}.
\end{lmm}

\subsection{Proof of Theorem \ref{T1}}
We now show how the ingredients above
(of which, Proposition \ref{pGZ} and Proposition \ref{Pkey} remain to be proved
later) are enough to prove Theorem \ref{T1}.

Let $\bp = \bp_n$ satisfy \re{regime} with limit $\btheta\in \Thfinite$.
Fix $J$ and take independent $U_1,\ldots,U_J$
with uniform $(0,1)$ distribution.
Proposition \ref{pGZ}
implies that as $n \to \infty$
\[ \sigma^{-1}(\bp)
 \left( G^{\bp}_I(U_{(1)}), \inf_{U_{(1)}\leq s \leq U_{(2)}} G^{\bp}_I(s),
  G^{\bp}_I(U_{(2)}), \inf_{U_{(2)}\leq s \leq U_{(3)}} G^{\bp}_I(s),
\ldots, G^{\bp}_I(U_{(J)}) \right)
\]
\[ \cd
 \left( Y^{\thv}(U_{(1)}), \inf_{U_{(1)}\leq s \leq U_{(2)}} Y^{\thv}(s),
  Y^{\thv}(U_{(2)}), \inf_{U_{(2)}\leq s \leq U_{(3)}} Y^{\thv}(s),
\ldots, Y^{\thv}(U_{(J)}) \right)
 . \]
By making the particular choice of $(\bp_n)$
used in Proposition \ref{Pkey},
\[ \sfrac{1}{2} \theta_0^2 \sp
 \left( H^{\bp}(U_{(1)}), \inf_{U_{(1)}\leq s \leq U_{(2)}} H^{\bp}(s),
  H^{\bp}(U_{(2)}), \inf_{U_{(2)}\leq s \leq U_{(3)}} H^{\bp}(s),
\ldots, H^{\bp}(U_{(J)}) \right)
\]
\begin{equation}
 \cd
 \left( Y^{\thv}(U_{(1)}), \inf_{U_{(1)}\leq s \leq U_{(2)}} Y^{\thv}(s),
  Y^{\thv}(U_{(2)}), \inf_{U_{(2)}\leq s \leq U_{(3)}} Y^{\thv}(s),
\ldots, Y^{\thv}(U_{(J)}) \right)
 . \lb{pHu}
\end{equation}
Appealing to Lemma \ref{Lrelate}, this implies
\[
\sfrac{1}{2}\theta_0^2 \sp \otimes {\cal S}^{\bp}_J
\cd
\TT^Y_J \]
where the right side
denotes the tree-with-edge-lengths obtained from sampling the
function $Y^{\thv}$,
and where convergence is the natural notion of convergence of
shapes and edge-lengths (\cite{me98} sec. 2.1).
Rescaling by a constant factor,
\[ \sp \otimes {\cal S}^{\bp}_J
\cd
\TT^{2 \theta_0^{-2}Y}_J . \]
But Proposition \ref{PwcICRT} showed
\[ \sp \otimes {\cal S}^{\bp}_J
\cd \TT^{\thv}_J \]
where the right side is the random tree-with-edge-lengths obtained
by sampling the ICRT $\TT^{\thv}$.
So we have established \re{TTZJ} and thereby proved Theorem \ref{T1} in the
case $\btheta\in \Thfinite$.

In the case $\sum_i\theta_i<\infty$, write $\btheta^n$ for the truncated
sequence $(\theta_0,\ldots,\theta_n,0,\ldots)$, and recall from
Lemma \ref{unifconv} that $Y^n=Y^{\sth^n}$ converges uniformly to $Y^{\sth}$.
By previous considerations this entails
$$\TT^{2\theta_0^{-2}Y^n}_J \cd \TT^{\sth}_J $$
for every $J\geq 1$. On the other hand, we have proved that
the left-hand term has the same law as
$c(\btheta^n)^{-1}\otimes\TT^{c(\sth^n)\sth^n}_J$ where
$c(\btheta^n)=(\sum_{0\leq i\leq n}\theta_i^2)^{-1/2}$ is the renormalization
constant so that $c(\btheta^n)\btheta^n\in \bTheta$. It thus remain to show
that this converges to $\TT^{\sth}$. Plainly the term $c(\btheta^n)$ converges
to $1$ and is unimportant. The result is then straightforward from the
line-breaking construction of the ICRT: $\TT^{\sth}_J$
can  be build out  of the  first (at  most) $2J$  points (cutpoints  and their
respective joinpoints) of the superimposition of 
infinitely many Poisson point processes on the line $(0,\infty)$. It is
easily checked that taking only the superimposition of the
$n$ first Poisson processes allows us to construct jointly
a reduced tree with same law as $c(\btheta^n)^{-1}\TT^{c(\sth^n)\sth^n}_J$ 
on the same probability space.
So for $n$ large the first $2J$ points of both point processes coincide and
we have actually $c(\btheta^n)^{-1}\TT^{c(\sth^n)\sth^n}_J=\TT^{\sth}_J$ on
this probability space. \cq

\rem Theorem \ref{T1} essentially consists of an ``identify the
limit'' problem, and that is why we are free to choose the
approximating $\bp_n$ in 
Proposition \ref{Pkey}.
But having proved Theorem \ref{T1}, we can reverse the proof above
to show that \re{pHu} holds true for {\em any} $\bp$ satisfying
\re{regime} with limiting $\btheta\in\Thfinite$. Indeed, the convergence
in \re{pHu} is equivalent to that of $\sp\otimes{\cal S}^{\bp}_J$ to
$\TT^{\sth}_J$ for every $J$.

\subsection{Skorokhod convergence of the discrete exploration process}

Suppose again that the ranked probability $\bp$ satisfies \re{regime}
with limit $\btheta\in \Thfinite$ with length $I$.
As observed in \cite{me98} (Theorem 5 and Proposition 7),
the convergence in \re{pHu} is equivalent to weak convergence of the
rescaled exploration process
to $Y^{\sth}$,
but using a certain topology on function space
which is weaker than the usual Skorokhod topology.
As noted in \cite{me98} Example 28,
assumption \re{regime}
is paradoxically not sufficient to ensure convergence
in the usual Skorokhod topology;
the obstacle in that example was the presence of exponentially
many (in terms of $1/\sp$) exponentially small $p$-values.
In this section we present some crude sufficient conditions
(\ref{H2},\ref{H1}); Proposition \ref{Pkey} will be a natural consequence
of the proof in section \ref{skocon}. The hypotheses are as follows.

First, we prevent very small $p$-values by making the assumption
\begin{equation}
1/p_* = o(\exp(\alpha/\sp))
\mbox{ for all } \alpha > 0
\lb{H2}
\end{equation}
where
\[ p_* := \min_i p_i . \]

Second, we will assume that most of the small
$p(\cdot)$-weights, as compared with the $I$ first, are of
order $\sp^2$.
Write $\bar{\bp}=(0,0,\ldots,p_{I+1},\ldots,p_n)$
for the sequence obtained from $\bp$ by truncating the first
$I$ terms. Let $\xi$ have distribution $\bp$ on $[n]$, and write
$\bar{p}(\xi )$ for the r.v.\ $\bar{p}_{\xi}$. We assume that there exists some
r.v.\ $0\leq Q<\infty$ such that the following ``moment generating function''
convergence holds:
\begin{equation}
\lim_{n\to \infty}E\left[\exp(\sfrac{\lambda\bar{p}(\xi)}{\sp^2}) \right]
=E\left[\exp\lambda Q\right]<\infty,
\lb{H1}
\end{equation}
for every $\lambda$ in some neighborhood of $0$. This implies that
$\bar{p}(\xi)/\sp^2\cd Q$, and also that the moments of all order exist
and converge to those of $Q$.

Then we have

\begin{thm}\label{usconv}
Suppose $\bp$ satisfies (\ref{regime}) with limit $\btheta\in\Thfinite$. Under
extra hypotheses (\ref{H2},\ref{H1}),
\begin{equation}\lb{H2B}
\sp H^{\bp}\cd \frac{2}{\theta_0^2}Y^{\sth}
\end{equation}
in the usual Skorokhod topology.
\end{thm}

\rem
The proof (section \ref{skocon}) rests upon applying the elementary large deviation
inequality
$P(S>s) \leq e^{-\lambda s} E \exp(\lambda S)$
to the independent sums involved in (\ref{Pt0},\ref{YpU}).
Hypothesis (\ref{H1}) is designed to make the application
very easy; it could surely be replaced by much weaker assumptions, such as
plain moment convergence conditions.

We would also guess that the convergence in \re{H2B} also holds with
$H^{\bp}$ replaced by
more general exploration processes, and in
particular the ``classical'' one, where each vertex $v$
is visited during an interval of length $1/n$ instead of $p_v$, or the
Harris (or contour) walk on the tree (see e.g.\ \cite[Chapter 2]{DGall02}).
We can easily verify the first
guess. Consider the $\bp$-tree $\psde(X_1,\ldots,X_n)$ defined as
in section \ref{subdfc} out of uniformly distributed independent r.v. Write
$w_1,\ldots,w_n$ for the vertices in depth-first order, and let $H^n(t)$ be
the height of the $w_i$ for which $i/n\leq t<(i+1)/n$ (and with the convention
$H^n(1)=H^n(1-)$).

\begin{crl}\label{clhp}
Suppose $\bp$ satisfies (\ref{regime}) with limit $\btheta\in\Thfinite$. Under
extra hypotheses (\ref{H2},\ref{H1}),
\begin{equation}\lb{H2Bn}
\sp H^n\cd \frac{2}{\theta_0^2}Y^{\sth}
\end{equation}
in the usual Skorokhod topology.
\end{crl}

\proof
By the functional weak law of large numbers for sampling without replacement,
we know that if $\pi$ is a uniform random permutation of the $n$ first
integers, the fact that
$\max_i p_{n,i}\to 0$ as $n\to \infty$ implies that if
$(S^0_n(t),0\leq t\leq 1)$ is the linear interpolation between points
$((i/n,\sum_{1\leq k\leq i}p_{\pi(i)}), 0\leq i\leq n)$ then $\sup_{0\leq
t\leq 1}|S^0_n(t)-t|\to0$ in probability. Now by the remark at the end 
of Sect.\ \ref{subdfc}, the cyclic order on vertices associated to the 
depth-first order
is uniform, so with the above notation for $i=w_1,\ldots,w_n$ the linear
interpolation $S_n$ between points
$((i/n,\sum_{1\leq k\leq  i} p_{w_k}),0\leq i\leq n)$ converges uniformly to
the identity in probability, since it is a (random) cyclic permutation of a
function distributed as
$S^0_n$. Noticing that $H^n=H^{\bp}\circ S^n$, the result follows. \cq

The convergence of the Harris walk follows from this
proposition by the arguments in \cite[Chapter 2.4]{DGall02}.



\section{Height profile}

This section is devoted to the proof of Theorem \ref{T2}. In this section, 
we do not assume that $\btheta\in\bTheta$ has finite length
nor that $\theta_0>0$.

\subsection{Continuity of the cumulative height profile}

We first prove the following intermediate lemma. Recall that
the {\em cumulative height process} of the
$\TT^{\sth}$ is defined as $\bar{W}^{\sth}(.)=
\mu^{\sth}\{v\in {\cal T}^{\theta}\, :\, \mbox{ht}(v)\leq.\}$,
where $\mu^{\sth}$ is the mass measure of $\TT^{\sth}$.

\begin{lmm}\label{continuity}
The cumulative
height process $\bar{W}^{\sth}$ is continuous for a.a.\ realizations of
$\TT^{\sth}$. Moreover, it has no flat interval, except its (possibly
empty) final constancy
interval, equal to  $[\sup_{v\in \TT^{\sth}}\HT(v),\infty)$.
\end{lmm}

\noindent{\em Proof of Lemma \ref{continuity}. }
Recall the recursive 
line-breaking construction of $\TT^{\sth}$ in the introduction, and
the  fact  from  \cite{me87}  that  the  tree  constructed  at  stage  $J$  is
distributed as the reduced tree $\TT^{\sth}_J$ of Sect.\ \ref{con-ptrees}. From this, we 
see that the leaves labelled $1,2,\ldots$ are a.s.\
at pairwise different heights, meaning that the measure $\d\bar{W}^{\sth}$
has no atom. Moreover, if $\bar{W}^{\sth}$ had a flat
interval (other than the final
constancy interval), this would mean that for some $h<\sup_{v\in\TT^{\sth}}$,
no leaf picked according to the mass measure can have a height in say
$(h-\epsilon,h+\epsilon)$ for some $\epsilon>0$.
But let $v$ be a vertex of
$\TT^{\sth}$ at height $h$. By the line-breaking construction, the fact that
branches have size going to $0$ and
the ``dense'' property of joinpoints, we can
find a joinpoint
$\eta^{*}$ at a distance
$<\epsilon/2$ of $v$ and so that the corresponding branch has length
$\eta<\epsilon/2$.
Since the leaves that are at the right-end of branches of the
line-breaking construction are distributed as independent sampled leaves from
the mass measure, this contradicts the above statement.
\cq

\subsection{Proof of Theorem \ref{T2}}

The reader can
consult \cite{kersting98m} for a similar treatment of convergence
of the height profile of Galton-Watson trees to a time-changed excursion
of a stable L\'evy process.

Suppose that $\bp={\bf p}^n$ satisfies the asymptotic regime \re{regime}.
Let ${\cal T}^{{\bf p}}$ be the ${\bf p}$-tree, and ${\cal T}^{{\sth}}$ the
limiting ICRT. Define $\bar{W}^{\sth}$ as above and recall the notation
$u(h)$ in \re{F1}. For $h\geq 0$ let
$$W^{\bp}(h)=\sum_{v\in {\cal T}^{{\bf p}},\,\HT(v)=
[\frac{h}{\sigma({\bf p})}]}
p_v=u\left(\left[\frac{h}{\sp}\right]+1\right)
-u\left(\left[\frac{h}{\sp}\right]\right),
\;\;\;\; h\geq 0$$
and $\bar{W}^{\bp}(h)=u([h/\sigma({\bf p})])$. Now let $U_1,U_2,\ldots$
be independent uniform$(0,1)$ random
variables. The sequence $((\bar{W}^{\bp})^{-1}(U_j),j\geq 1)$
has the law of the
heights of an i.i.d.
random sample of vertices of ${\cal T}^{{\bf p}}$, chosen according to $\bp$,
and the same holds for
$((\bar{W}^{\sth})^{-1}(U_j),j\geq1)$ and the tree
${\cal T}^{\theta}$, with the mass measure $\mu^{\sth}$
as common law. For $J\geq 1$ let $\bar{W}^{\bp}_J(h)$
be the associated empirical distribution of the first $J$ terms, defined by
$$\bar{W}^{\bp}_J(h)=\frac{1}{J}\sum_{i=1}^{J}\ind_{\{(\bar{W}^{\bp})^{-1}
(U_i) \leq
h\}}=\frac{1}{J}\sum_{i=1}^{J}\ind_{\{U_i\leq \bar{W}^{\bp}(h)\}},$$
and define $\bar{W}^{\sth}_J(h)$ in a similar way.

By Proposition \ref{PwcICRT},
we have that the random Stieltjes measure
$\d\bar{W}^{\bp}_J$ converges in law to $\d\bar{W}^{\sth}_J$ as $n\to\infty$ for
every $J\geq 1$.
Moreover, the empirical measure of an i.i.d.\ $J$-sample of leaves
distributed according to $\mu^{\sth}$ converges to $\mu^{\sth}$, implying
$\d\bar{W}^{\sth}_J\cd \d\bar{W}^{\sth}$ as $J\to \infty$. Thus,
for $h\geq 0$ and $J_n\to\infty$ slowly enough,
$$\d\bar{W}^{\bp}_{J_n}\build\to_{}^{d}\d\bar{W}^{\sth}.$$

Now let $F_J(x)=J^{-1}\sum_{i=1}^{J}\ind_{\{U_i\leq x\}}$ be the
empirical distribution associated
to the uniform variables $U_1,\ldots,U_J$.
Then $\sup_{h\geq
0}|\bar{W}^{\bp}_J(h)-\bar{W}^{\bp}(h)|\leq\sup_{x\in[0,1]}|F_J(x)-x|$,
which by the Glivenko-Cantelli
Theorem converges to $0$ as $J\to\infty$, and this convergence is
uniform in $n$. Hence
the random measure $\d\bar{W}^{\bp}$ converges in distribution
to $\d\bar{W}^{\sth}$ for the weak topology on measures.
Thanks to Lemma \ref{continuity} we may improve this to
$$\bar{W}^{\bp}(\cdot)\build\to_{}^{d}\bar{W}^{\sth}(\cdot)$$
where the convergence is weak convergence of processes for the topology
of uniform convergence.
It is then an elementary consequence of Lemma \ref{continuity}
that $\bar{W}^{\bp}((\bar{W}^{\bp})^{-1}(\cdot))$ converges in law for the 
uniform convergence topology to the identity function on $[0,1]$.

Equation \re{F2} can be rewritten as
\begin{equation}\label{F2r}
W^{\bp}(h)=
\fexp(\bar{W}^{\bp}(h)), \hskip1cm h\geq 0,
\end{equation}
so the convergence in distribution of $\bar{W}^{\bp}$, the fact that
its limit is strictly increasing and continuous, and \re{pFX} imply
that the sequence of random processes $(\sp^{-1}W^{\bp})$ is tight. Thus,
the pair $(\sp^{-1}W^{\bp},\bar{W}^{\bp})$ is tight, and
up to extraction of a subsequence, we can
suppose that $(\sp^{-1}W^{\bp},\bar{W}^{\bp})\cd (W,\bar{W}')$
for some process $W$, and where $\bar{W}'$ has the same law as
$\bar{W}^{\sth}$. Suppose further
by Skorokhod's embedding theorem that the convergence
is almost-sure. By definition
$$\int_0^h\frac{W^{\bp}(u)}{\sp}\d u=\bar{W}^{\bp}(h-\sp)+R(n,h)$$
where $R(n,h)\leq\bar{W}^{\bp}(h)-\bar{W}^{\bp}(h-\sp)$
goes to $0$ uniformly as $n\to\infty$ by
continuity of the limiting $\bar{W}'$. So necessarily,
$$\int_0^hW(u)\d u=\bar{W}'_h, \hskip1cm h\geq 0$$
for every $h\geq 0$, so that the only possible limit
$W$ is the density of $\d\bar{W}'$.
Therefore, the height profile $W^{\sth}$ of
the ICRT exists and $(\sp^{-1}W^{\bp},\bar{W}^{\bp})\cd(W^{\sth},
\bar{W}^{\sth})$.
Looking back at \re{F2r} we have
$$\sp^{-1}W^{\bp}((\bar{W}^{\bp})^{-1}(u))=\sp^{-1}
\fexp(\bar{W}^{\bp}((\bar{W}^{\bp})^{-1}(u))), \hskip1cm 0\leq u\leq 1,$$
so by the convergence of $\bar{W}^{\bp}((\bar{W}^{\bp})^{-1}(\cdot))$ and
\re{pFX}, we obtain convergence in distribution of
the right-hand side to
$X^{\thv}$. By the convergence in law of $W^{\bp}$
this finally implies that
$W^{\sth}((\bar{W}^{\sth})^{-1}(\cdot))\ed
X^{\thv}(\cdot)$ and Theorem \ref{T2} is proved. \cq


\noindent{\em Proof of Corollary \ref{httree}. }
By the proof of Lemma \ref{continuity}, the only constant interval of the
width process of the ICRT is $[\sup_{v\in\TT^{\sth}}\HT(v),\infty)$. Thus
the height of the tree, $\sup_{v\in\TT^{\sth}}\HT(v)$, is the first point
after which the width process remains constant. By (\ref{Lamperti}), 
this point has same law as $\int_0^1 \d s/X^{\sth}_s$. \cq


\section{The exploration process}\label{htprocess}

To shorten notation, for $A \subseteq [n]$ we write
$p(A)$ for 
the quantity $\sum_{j\in A}p_j$.

\subsection{Proof of Proposition 2}
\label{convY}

Let $\bp$ satisfy \re{regime} for some limiting $\btheta\in \Thfinite$, with
length $I$. In this subsection we suppose that the $\bp$-tree $\tree^{\bp}$
is constructed from the
process $\fexp$ by the depth-first search construction of section
\ref{con-ptrees}.
Moreover, since we have \re{pFX}
the convergence in law $\sp^{-1}\fexp\cd X^{\sth}$,
we suppose by Skorokhod's representation theorem that our probability
space is such that the convergence holds almost surely.
Recall that in the depth-first search construction of the $\bp$-tree out
of the process $\fexp$, the $i$-th examined vertex $v=w_i$ is examined
during an interval $[e(v)-p_v,e(v))$, during which the labels of jumps of
$\fexp$ determine the set ${\cal B}_v$ of children of $v$.

We begin with two useful observations. First, if $v$ is a vertex
of $\TT^{\bp}$ and
if $\TT^{\bp}_v$ denotes the {\em fringe subtree} of $\TT^{\bp}$ rooted at 
$v$, that is, the subtree of descendents of $v$, then for every vertex $w$ of
$\TT^{\bp}_v$ one has
\begin{equation}\lb{minw}
\fexp(e(w))\geq \fexp(e(v))-p({\cal B}_v).
\end{equation}
To argue this, simply recall formula \re{deprop} and notice that
${\cal N}(v)\subseteq{\cal N}(w)\cup{\cal B}_v$.

Second, notice that since $\max_j p_j\to0$ and
the limiting process $X^{\sth}$ is continuous
except for a finite number $I$ of upward jumps, we must necessarily
have that a.s.\ as $n\to\infty$,
\begin{equation}\lb{domd}
\eta_n:=\max_{v\in[n]}\left|\inf_{u\in[e(v)-p_v,e(v))}
(\fexp(u)-\fexp(e(v)-p_v))\right| =o(\sp).
\end{equation}

\begin{lmm}\label{limB}
Almost surely
$$\max_{j\in[n]}\sp^{-1}\left|p_j-p({\cal B}_j\setminus[I])\right|\to 0.$$
\end{lmm}

\proof
As mentioned, for every vertex $v\in[n]$,
$$\fexp(e(v))-\fexp(e(v)-p_v)=p({\cal B}_v)-p_v.$$
Consider the process $F^{\bp\downarrow}$ defined by
$$F^{\bp\downarrow}(s)=\fexp(s)-\sum_{1\leq i\leq I}p_i\ind_{\{s\geq x_i'\}}$$
where as above $x'_i$ is the time when $\fexp$ has its jump with size $p_i$.
Easily, $\sp^{-1}F^{\bp\downarrow}$ converges in the Skorokhod space to
the process $X^{\sth\downarrow}$ defined by
$$X^{\sth\downarrow}_s=X^{\sth}_s-\sum_{1\leq i\leq I}\theta_i\ind_{\{s\geq t_i\}}$$
where $t_i$ is the time when $X^{\sth}$ jumps by $\theta_i$. This process is
continuous, hence $\max_j p_j\to0$ implies
$$\sp^{-1} \max_v |F^{\bp\downarrow}(e(v))-F^{\bp\downarrow}(e(v)-p_v)|\to 0.$$
Now the quantity $F^{\bp\downarrow}(e(v))-F^{\bp\downarrow}(e(v)-p_v)$ equals
\begin{eqnarray*}
&&\fexp(e(v))-\fexp(e(v)-p_v)-\sum_{1\leq i\leq I}p_i
\ind_{\{x'_i\in(e(v)-p_v,e(v)]\}} \\
&=&p({\cal B}_v)-p_v-p({\cal B}_v\cap[I])
\end{eqnarray*}
implying the lemma. \cq

Now, for $v$ a non-root vertex of $\TT^{\bp}$
let $f(v)$ be its parent. For $i\in[I]$ and
$n$ large enough, $i$ is not the root (since the limiting $X^{\sth}$ does
not begin with a jump), so $f(i)$ exists.

\begin{lmm}\label{mi}
Let $i \in I$.
Let ${\cal M}(i)$ be the set of
descendents of $f(i)$ that come strictly before $i$ in depth-first order.
Suppose that $f(i)\notin[I]$ for $n$ large enough.
Then as $n\to \infty$,  $p({\cal M}(i))\to0$ almost surely.
\end{lmm}

\proof
A variation
of \re{minw} implies for any
$v\in{\cal M}(i)$ and $n$ large that
\begin{equation}\lb{interm}
\fexp(e(v))\geq \fexp(e(f(i)))-p({\cal B}_{f(i)}\setminus[I]).
\end{equation}
Indeed, it is clear that for $n$ large the sets ${\cal B}_v\cap[I]$ contain
at most one element, otherwise the
Skorokhod convergence $\sp^{-1}\fexp\to X^{\sth}$ would fail as two or
more upward jumps of non-negligible sizes could occur in an ultimately
negligible interval. Moreover, for $v\in{\cal M}(i)$, it is clear that
${\cal N}(v)$ contains $i$, hence \re{interm}. Thus
$$\inf_{e(f(i))\leq u\leq e(f(i))+p({\cal M}(i))}\fexp(u)
\geq \fexp(e(f(i)))-p({\cal B}_{f(i)}\setminus[I])-\eta_n,$$
with $\eta_n$ defined at \re{domd},
since the vertices of ${\cal M}(i)$ are visited during the
interval $[e(f(i)),e(f(i))+p({\cal M}(i))]$. Since $\sp^{-1}
\fexp(e(f(i)))$ is easily
seen to converge to $X^{\sth}_{t_i}$, by \re{domd}, Lemma \ref{limB} and
the fact that $f(i)\notin[I]$ for $n$ large,
if $p({\cal M}(i))$ did not converge to
$0$, by extracting along a
subsequence we could find an interval $[t_i,t_i+\eps]$
with $\eps>0$ where $X^{\sth}_u\geq X^{\sth}_{t_i}$, and this is a.s.\
impossible by Lemma \ref{XandY}. \cq

The assumption that $f(i)\notin[I]$ may look strange since it is intuitive 
that the child of some $i\in[I]$ is very unlikely to be in $[I]$ for $n$ large
(e.g.\ by Theorem \ref{T2}). We actually have:

\begin{lmm}\label{limtheta}
For every $i\in[I]$, almost surely, ${\cal B}_i\cap[I]=\emptyset$
for $n$ large, and
$$\sp^{-1}p({\cal B}_i)\to \theta_i.$$
\end{lmm}

\proof
By Lemma \ref{limB} it suffices to prove that a.s.\ for large $n$,
${\cal B}_i\cap[I]=\emptyset$. Suppose that there exist $i,j\in[I]$ such
that $j$ is the child of $i$ in the $\bp$-tree infinitely often. Since
$I<\infty$, we may further suppose that $f(i)\notin[I]$ by taking
(up to extraction) the least
such $i$ in depth-first order.
By definition, $\fexp$ has
a jump with size $i$ in the interval $[e(f(i))-p_{f(i)},e(f(i))]$.
Moreover, it follows from the definition of ${\cal M}(i)$ that
$e(i)-p_i=e(f(i))+p({\cal M}(i))$. Since
the vertex $i$ is examined in the interval
$[e(i)-p_i,e(i)]$ and $p({\cal M}(i))\to 0$ by the preceding lemma, the fact
that $f(j)=i$ implies that the
jumps with
size $p_i$ and $p_j$ occur within a vanishing interval
$[e(f(i))-p_{f(i)},e(i)]$. Therefore, the
Skorokhod convergence of $\sp^{-1}\fexp$ to $X^{\sth}$ would fail. \cq

Now recall the definition (\ref{ri}) of the processes $r_i^{\bp}$ used to build $\gip$
in section \ref{ptree-icrt}, and that $x'_i$ is the time when $\fexp$ jumps
by $p_i$. .

\begin{lmm}\label{convri}
For every $i\in[I]$, as $n\to \infty$, we have
$$\sp^{-1}\left|
\inf_{x'_i\leq u\leq s}\fexp(u)-\fexp(x'_i-)-r_i^{\bp}(s)\right|\to 0$$
a.s.\ uniformly in $s\in[x_i',e(i)+p(\TT^{\bp}_i)]$.
\end{lmm}

\proof
Let $i\in[I]$, and let ${\cal B}_i=\{v_1,v_2,\ldots,v_k\}$
(with $k=|{\cal B}_i|$) where
$v_1,v_2,\ldots$ are in depth-first order. For $1\leq j\leq k$ let
also $v_j'$ be the last examined vertex of $\TT^{\bp}_{v_j}$ in depth-first
order, that is, the predecessor of $v_{j+1}$ if $j<k$. Then one has, for
every $1\leq j\leq k$ and $w\in \TT^{\bp}_{v_j}$
$$\fexp(e(w))\geq \fexp(e(v_j))-p({\cal B}_{v_j}),$$
as follows from \re{minw}. Rewrite this as
$$\fexp(e(w))\geq \fexp(e(i))-\sum_{1\leq r\leq j-1}p_{v_r}$$
and check that the right hand side equals
$\fexp(e(v_{j-1}'))$.
In particular, we obtain
$$\left|\inf_{v:e(v)\in[e(i),e(w)]}\fexp(e(v))
-\fexp(e(i))+\sum_{1\leq r\leq j-1}p_{v_r}\right|
\leq \max_{1\leq j\leq k}p_{v_j} . $$
Now check that for $w$ a vertex of $\TT^{\bp}_{v_j}$, one has
$r_i^{\bp}(e(w))=\sum_{j\leq r\leq k}p_{v_r}$. For $s$ as in the
statement of the lemma deduce, for $n$ large (since ${\cal B}_i\cap[I]=
\emptyset$ by Lemma \ref{limtheta}),
$$\left|\inf_{u\in[x'_i,s]}\fexp(u)-\fexp(e(i))+p({\cal B}_i)-r_i^{\bp}(u)
\right|\leq 2\max_{j\notin[I]}p_{j}+\eta_n+\eta_n'$$
where
$$\eta_n'=\max_{x'_i\leq u\leq e(i)}|\fexp(u)-\fexp(e(i))|$$
which
is $o(\sp)$ by Lemma \ref{mi} and the convergence $\sp^{-1}\fexp\to X^{\sth}$.
We conclude, using the fact that $\sp^{-1}\fexp(x_i'-)\to X^{\sth}_{t_i-}$,
which is equal to the limit of $\sp^{-1}(\fexp(e(i))-p({\cal B}_i))$,
as follows from Lemmas \ref{mi} and \ref{limtheta}. \cq

\noindent{\em Proof of Proposition \ref{pGZ}. }
We prove that the process $\sp^{-1}r_i^{\bp}$
converges to
the $R^{\sth}_i$ of section \ref{sec-cXZ} in the Skorokhod topology, for every
$i$. In view  of Lemma \ref{convri}, and since by definition  of $r_i$ one has
$r_i(u)=0$ for $u\geq e(i)+p(\TT^{\bp}_i)$, the only thing to do is to show 
that $e(v'_k)=e(i)+p(\TT^{\bp}_i)$
converges to the
$T_i$ of section \ref{sec-cXZ}. Since 
$e(v'_k)\geq  \inf\{s\geq  x'_i:r^{\bp}_i(s)=0\}$,  we  obtain  that  $\liminf
e(v'_k)\geq T_i$. Suppose $\ell=\limsup e(v'_k)>T_i$, and up to extraction 
suppose that $\ell$ is actually the limit of $e(v'_k)$.  From the fact that 
$\fexp(e(v'_k))=\fexp(e(i))-p({\cal B}_i)$, hence $\sp^{-1}\fexp(e(v'_k))$
converges to $X^{\sth}_{t_i-}$ by Lemmas \ref{mi} and \ref{limtheta}, we
would find $\ell>T_i$ with 
$X^{\sth}_{\ell}=X^{\sth}_{t_i-}$ and $X^{\sth}_s\geq X^{\sth}_{t_i-}$
for $s\in[T_i,\ell]$, and this is almost surely
impossible by Lemma \ref{XandY} as $X^{\sth}_{t_i-}$ would be a local minimum 
of $X^{\sth}$, attained at time $T_i$. 

Without extra argument we cannot conclude that the sum
$\sp^{-1}(\fexp-\sum_{i=1}^Ir_i^{\bp})$ converges to $X^{\sth}-\sum_{i=1}^I
R^{\sth}_i$, but this is nonetheless
true for the following reason.
The process $R^{\sth}_i$ is continuous except for
one jump at $t_i$, and the process $r_i^{\bp}$ has precisely one jump with 
size
$p({\cal B}_i)$ at time $x'_i$, that is, at the {\em same} time as the jump of
$\fexp$ with size $p_i$. Together with Lemma \ref{limtheta}, we obtain the 
Skorokhod convergence $\sp^{-1}\gip\to Y^{\sth}$. \cq

\subsection{Proof of Theorem \ref{usconv}}\label{skocon}

As above, we suppose that $\bp$ is a ranked probability distribution
satisfying \re{regime} for some limiting $\btheta$ with length $I$, and
we suppose that the $\bp$-tree $\TT^{\bp}$ is obtained by the depth-first
construction of section \ref{con-ptrees} out of the process $\fexp$.
We are going to show the following result:
\begin{prp}\label{P1}
Under extra hypotheses
(\ref{H2},\ref{H1})
on $(\bp^{(n)})$,
as $n \to \infty$
\[ \max_v
\left| \sfrac{\theta_0^2\sp}{2} \HT(v)
- \sp^{-1} \gip(e(v))
\right|
\cp 0 . \]
\end{prp}
We first show how Theorem \ref{usconv} and Proposition \ref{Pkey}
are easy consequences of Proposition \ref{P1}.

\noindent{\em Proof of Theorem \ref{usconv}. }
Since $\sp^{-1}\gip$ converges uniformly in distribution
to a continuous process, and since $H^{\bp}$ does not vary in the
intervals $[e(v)-p_v,e(v))$, the last displayed convergence extends
to
\begin{equation}\lb{smg}
\max_{u\in[0,1]}\left|\sfrac{\theta_0^2\sp}{2}H^{\bp}(u)-\sp^{-1}\gip(u)
\right|\cp 0,
\end{equation}
and then Proposition \ref{pGZ}
implies Theorem \ref{usconv}. \cq

\noindent{\em Proof of Proposition \ref{Pkey}. }
For Proposition \ref{Pkey}, we choose the
following approximating sequence $\bp^{(n+I)}$ for $\btheta\in\Thfinite$
with length $I$. Given $n$, let
$z_n=\sqrt{n}/\theta_0$, $s_n=n+z_n\sum_{1\leq i\leq I}\theta_i$ and
\begin{equation}\lb{particular}
\left\{ \begin{array}{lr}
p_i=\frac{z_n\theta_i}{s_n} & \mbox{ if }1\leq i\leq I\\
p_i=\frac{1}{s_n} & \mbox{ if }I+1\leq i\leq n+I .
\end{array} \right.
\end{equation}
It is trivial to see that this sequence fulfills hypotheses
(\ref{H2},\ref{H1}). Hence \re{smg} is satisfied, and Proposition \ref{Pkey}
is an immediate consequence. \cq

We now mention three consequences of hypotheses (\ref{H2},\ref{H1}) that
will be used later. First, notice that $p_*\leq 1/n$ since
$\bp$ is a probability on $[n]$, so \re{H2} implies
\begin{equation}
n=o(\exp(\alpha/\sp)) \mbox{ for all }\alpha>0.
\lb{H2bis}
\end{equation}
Second, \re{H1} implies convergence of all moments of $\bar{p}(\xi)/\sp^2$,
and in particular
\begin{eqnarray}\lb{expecQ}
E\left(\frac{\bar{p}(\xi)}{\sp^2}\right)
&=&\sum_{i\notin[I]}p_i^2/\sp^2 \nonumber\\
&=&1-\sum_{i\in[I]}p_i^2/\sp^2\nonumber\\
&\build\to_{n\to \infty}^{}&\theta_0^2=E(Q).
\end{eqnarray}
Third, for every $\lambda$ in a neighborhood of $0$,
\begin{equation}
\sp^2 \sum_{i\notin[I]} \left[ \exp\left(\sfrac{\lambda p_i}{\sigma^2}
\right) - 1 - \sfrac{\lambda p_i}{\sigma^2}  \right]
\build\to_{n\to\infty}^{} E \sfrac{1}{Q} \left[ \exp( \lambda Q
) - 1 - \lambda Q \right]
< \infty.
\lb{H4}
\end{equation}
Indeed, the left side can be rewritten as
$E\left(\frac{\sp^2}{ \bar{p}(\xi)} \left[ \exp(\frac{\lambda
\bar{p}(\xi)}{\sigma^2})
- 1 - \sfrac{\lambda \bar{p}(\xi)}{\sigma^2} \right]\right)$, where
the function
$f(x)=(e^{\lambda x}-1-\lambda x)/x$ is understood to equal its limit $0$ 
at $0$.
Since it is bounded in a neighborhood of $0$ and
dominated by $e^{\lambda x}$ near $\infty$, the convergence of this expectation
is an easy consequence of \re{H1}.

The first step in the proof of Proposition \ref{P1}
is to relate $H(\cdot)$ to another function 
${\cal G}(\cdot)$ measuring
``sum of small $\bp$-values along path to root''. Let ${\cal A}(v)$ be
the  set of ancestors of $v$ in the $\bp$-tree, and let
\begin{equation}
{\cal G}(v) := p({\cal A}(v)\setminus[I]) .
\end{equation}
\begin{lmm}
\label{L1}
Under extra hypotheses (\ref{H2},\ref{H1}),
as $n \to \infty$
for fixed $K>0$
\[ \max_{v: \HT(v) \leq K/\sp} \left|
\sp \theta_0^2 \, \HT(v) - \sp^{-1}{\cal G}(v)
\right| \cp 0 . \]
\end{lmm}
\proof
Let $V$ be a $\bp$-distributed random vertex.
Fix $\eps > 0$.
It is enough to prove
that as $n \to \infty$
\[
P \left( |
\sp\theta_0^2\, \HT(V)  - \sp^{-1} {\cal G}(e(V))
| > \eps, \ \sp \HT(V) \leq K \right)
= o(p_*) . \]
Let $\xi$ have distribution $\bp$ on $[n]$ and let
$(\xi_i, i \geq 1)$ be i.i.d.
By the ``birthday tree'' construction of the $\bp$-tree
\cite[Corollary 3]{jpmc97b}
we have equality of joint distributions
\[ (\HT(V), {\cal G}(V)) \ed (T-2, \sum_{i=1}^{T-1} \bar{p}(\xi_i)) \]
where
\[ T:= \min \{j \geq 2: \xi_j = \xi_i \mbox{ for some } 1 \leq i < j \} \]
is the first repeat time in the sequence $\xi_i$.
So it is enough to prove
\[ P \left(
\left| \sp \theta_0^2 (T-2) - \sp^{-1} \sum_{i=1}^{T-1} \bar{p}(\xi_i) \right|
> \eps, \ \sp (T-2) \leq K \right) = o(p_*) . \]
We may replace $T-2$ by $T-1$ and $\theta_0^2$ by
$E(\sfrac{\bar{p}(\xi)}{\sp^2})$ by the above remark.
Rewriting in terms of
$\tilde{p}(i) := \frac{\bar{p}(i)}{\sp^2} - E(\frac{\bar{p}(\xi)}{\sp^2})$,
we need to prove
\[ P \left( \left|
\sum_{i=1}^{T-1} \tilde{p}(\xi_i) \right| > \eps/\sp
, \ T-1 \leq K/\sp
\right) = o(p_*) . \]
Now we are dealing with a mean-zero random walk, and
classical fluctuation inequalities
(e.g. \cite{durrett95} Exercise 1.8.9)
reduce the problem to proving the fixed-time bound
\begin{equation}
 P \left( \left|
\sum_{1 \leq i \leq K/\sp} \tilde{p}(\xi_i)
\right| \geq \eps/\sp \right)
= o(p_*) . \lb{Pt0}
\end{equation}
We now appeal to assumption (\ref{H1}),
which basically says that the sums in question behave as
if the summands had distribution
$Q -\theta_0^2$ not depending on $n$.
More precisely, the elementary large deviation inequality applied to
the probability in \re{Pt0} but without the absolute values implies that
for any small $\lambda>0$,
$$
\log  P \left(
\sum_{i=1}^{K/\sp} \tilde{p}(\xi_i)
\geq \eps/\sp \right)\leq   -\frac{\lambda\eps}{\sp}+\frac{K}{\sp}
\log(E(\exp(\lambda\tilde{p}(\xi))).$$
Assumption \re{H1} and the convergence of the expectation of $\bar{p}(\xi)$
allows us to rewrite the log term on the right as
$$\frac{K}{\sp}\log E(\exp(\lambda(Q-\theta_0^2)))+
\frac{K\eta_{\lambda}(n)}{\sp E(\exp(\lambda(Q-\theta_0^2)))},$$
where $\eta_{\lambda}(n)\to 0$ as $n\to \infty$ for any fixed $\lambda$.
We now choose $\lambda$ small enough so that $-\lambda\eps+
K\log E(\exp(\lambda(Q-\theta_0^2)))=-\delta<0$ and we let $n \to \infty$,
obtaining the bound $\exp(-\delta'/\sp)$, for some $\delta'>0$,
for the probability in \re{Pt0}
without absolute values, but the other side of the inequality is similar.
Now assumption \re{H2} gives the desired bound \re{Pt0}. \cq


The next, rather strange-looking lemma does most of the
work in relating the processes $\gip(\cdot)$ and ${\cal G}(\cdot)$.

Given a probability distribution $\bp$ on $[n]$ and given
a subset $A \subset [n]$, let
$\bq$ be the probability distribution obtained by lumping the
points $A$ into a single point; that is,
$q_1 = p(A)$
and the multiset $\{q_i, i \geq 2\}$
is the multiset $\{p_i, i \not\in A\}$.
We also let ${\cal I}$ be the set of ``large'' $\bq$-values, except $q_1$. 
Precisely, ${\cal I}$ is such that the multisets
$\{p_v,v\in[I]\setminus A\}=\{q_v,v\in{\cal I}\}$ are equal. Then

\begin{lmm}
\label{L2}
Suppose $\bp=\bp^{(n)}$ satisfies the regime
\re{regime} and extra hypotheses (\ref{H2},\ref{H1}).
Let $A = A^{(n)} \subset [n]$ and define $\bq$ as above.
Define a random variable $X = X(\bq)$ as follows.
Take a $\bq$-tree, condition on vertex $1$ being the root. Let ${\cal B}_1$
be the set of children of $1$, and
for each $v\in{\cal B}_1$ toss two coins $c_1$ and $c_2$,
$c_1$ a fair coin and $P(c_2=\mbox{Heads})=p(A\setminus[I])/p(A)$, and set
\[ X := \sum \{q_v: \ v\in{\cal B}_1\setminus{\cal I},
\mbox{ coins $c_1$ and $c_2$ land Heads} \} . \]
Suppose $q_1 \leq K\sp$ and set $\bar{q}_1=p(A\setminus[I])$. 
Then for fixed $\eps > 0$ there exists $\delta=\delta(\eps,K)>0$ with
\[ P( |X - \sfrac{1}{2}\bar{q}_1| > \eps \sp) \leq \exp(-\delta/\sp)=o(1/n), \]
where the $o(1/n)$ is thus uniform over $q_1\leq K\sp$.
\end{lmm}
\proof
Consider the random variable
\[ Y:= \sum_{i \not\in A\cup[I]} p_i \ind_{\{U_i \leq \bar{q}_1/2\}} \]
where the $(U_i)$ are independent uniform$(0,1)$.
The key relation is
\begin{equation}
P(X \in \cdot) \leq \sfrac{1}{q_1} P(Y \in \cdot) .
\lb{key-rel}
\end{equation}
This follows from the {\em breadth-first} construction of the $\bp$-trees.
In that construction of a $\bq$-tree, vertices $i$ are associated
with uniform$(0,1)$ r.v.'s $U^\prime_i$ in such a way that, if
vertex $1$ happens to be the root, then the children $v$ of $1$
are the vertices $v$ for which
$U_v := U^\prime_v - U^\prime_1 \bmod 1$ falls within $(0,q_1)$.
Thus, writing
\[ X' := \sum \{q_v: v \in{\cal B}_1\setminus{\cal I}\} \]
\[ Y':= \sum_{i\notin A\cup[I]} p_i \ind_{\{U_i \leq q_1\}} \]
we have
\[ X' = Y' \mbox{ on the event } \{
\mbox{ vertex $1$ is root } \} . \]
So
\[ P(X' \in \cdot| \mbox{ $1$ is root})
\leq \frac{P(Y' \in \cdot)}{P(\mbox{ $1$ is root})}
= \sfrac{1}{q_1} P(Y' \in \cdot) . \]
The stated inequality (\ref{key-rel}) follows by applying an independent
Bernoulli$(\bar{q}_1/(2q_1))$ thinning procedure to both sides.

Now write $c = \bar{q}_1/2$ and let us study the centered version of $Y$:
\begin{equation}
\tilde{Y} := \sum_{i \not\in A\cup[I]} p_i(\ind_{\{U_i \leq c\}} - c) .
\lb{YpU}
\end{equation}
The elementary large deviation bound, applied to
$\tilde{Y}/\sp^2$, is: for arbitrary
$\lambda > 0$,
\[
\log P(\tilde{Y} > \eps \sp)
\leq \frac{-\lambda \eps}{\sp}
+ \log E \exp(\lambda \tilde{Y}/\sp^2) .
\]
We calculate
\begin{eqnarray*}
\lefteqn{
 \log E \exp(\lambda \tilde{Y}/\sp^2) }\\
&=&
\sum_{i \not\in A\cup[I]} \left\{
\sfrac{-\lambda p_i}{\sp^2} c
+ \log \left[ 1 + c(
e^{\lambda p_i/\sp^2} -1)
\right] \right\} \\
&\leq& c
\sum_{i\in[n]} \left\{
e^{\lambda p_i/\sp^2} -1
- \sfrac{\lambda p_i}{\sp^2}
\right\},
\end{eqnarray*}
since the quantities we are summing are positive,
and by (\ref{H4}) the bound is asymptotic to
$c\sp^{-2}\Phi(\lambda)$ for
\[ \Phi(\lambda) :=
 E \sfrac{1}{Q} \left[ \exp(\lambda Q) - 1 - \lambda Q \right]
 . \]
By hypothesis $c := \bar{q}_1/2 \leq K\sp$, so $c\sp^{-2} \leq K\sp^{-1}$. 
So there is a constant $C_1=C_1(K)$ such that
\[
\log P(\tilde{Y} > \eps \sp)
\leq \sfrac{1}{\sp}
\left( - \lambda \eps + C_1 \Phi(\lambda) \right) . \]
But $\Phi^\prime(0) = 0$ and so
$\Phi(\lambda) = o(\lambda)$ as $\lambda \downarrow 0$,
so the right side is strictly negative for small $\lambda > 0$.
So there exists $\delta_1 = \delta_1(\eps,K) > 0$ such that
\[ P(\tilde{Y} > \eps \sp)  \leq \exp(-\delta_1/\sp) . \]
Since $Y - \tilde{Y} = c \sum_{i \not\in A\cup[I]} p_i \leq \bar{q}_1/2$
we have established the one-sided inequality
\[ P(Y - \sfrac{1}{2}\bar{q}_1 > \eps \sp) 
\leq \exp(-\delta_1/\sp) . \]
The other side of the inequality is similar
except for this last step: we cannot bound so easily the quantity
$\tilde{Y}-Y$. However, by \re{regime},
$$\sum_{i\notin A\cup[I]}p_i=
1-p(A\cup[I])\geq 1-q_1-\sum_{i=1}^Ip_i\geq 1-C_2\sp $$
for some
$C_2=C_2(K)<\infty$. Thus $Y-\tilde{Y}\geq c(1-C_2\sp)$ and we can conclude 
as above by the existence of $\delta_2=\delta_2(\eps,K)$ satisfying 
$$P(\sfrac{1}{2} \bar{q}_1-Y>\eps \sp)\leq \exp(-\delta_2/\sp).$$ 
So, letting $\delta'=\delta_1\wedge \delta_2$, 
\[ P(|Y - \sfrac{1}{2}\bar{q}_1 | > \eps \sp) \leq
2\exp(-\delta'/\sp) . \]
Now (\ref{key-rel}) and hypothesis (\ref{H2}) and its consequence \re{H2bis}
establish Lemma \ref{L2} (with any $\delta<\delta'$). \cq

For the next lemma, recall the definition of ${\cal N}(v)$ around \re{deprop}
and let
${\cal N}^*(v)$ be the subset of vertices of ${\cal N}(v)$ which are not in
$[I]$ and whose parent is not in $[I]$ either.
\begin{lmm}
\label{L3}
Fix $j \in [n]$ and a subset $A \subset [n]$ with $j \in A$.
Take a random $\bp$-tree and condition on ${\cal A}(j)=A$.  Let also
$v_1,\ldots,v_k$ be the children of $j$ that are not in $[I]$
and let $c^*(j)=\sum_{1\leq l\leq k}b_lp_{v_l}$, where the $b_l$'s are 
independent Bernoulli  random variables  with parameter $1/2$,  independent of
the $\bp$-tree. Define
\[ X^*:= p({\cal N}^*(j))-c^*(j) . \]
Then $X^*$ is distributed as the random variable $X$ in Lemma \ref{L2}.
\end{lmm}
\proof
Order $A$ as $v_0,v_1,\ldots, j$, arbitrarily except for ending with $j$.
Let $\bT^*$ be the set of rooted trees on $[n]$ with root $v_0$
whose path to $j$ is the path $v_0,v_1,\ldots,j$.
Let $\bT^\oplus$ be the set of rooted trees on
$[n] \setminus A \cup\{\oplus\}$
with root $\oplus$.
There is a natural map $\bT^* \to \bT^\oplus$:
``lump the vertices in $A$ together into a single vertex $\oplus$''.
It is straightforward to check,
from the combinatorial definition (see e.g.\ \cite{jp01hur}) of $\bp$-tree,
that this map takes
the distribution of $\bp$-tree (conditioned to $\bT^*$)
into the distribution of a $\bq$-tree (conditioned on having root
$\oplus$). Also, we have the extra constraint in $X^*$ that the parents of 
the vertices we are summing on are not in $[I]$, but conditionally
on the fact that $v$ has some parent in $A$, it is easy that the
parent is in $[I]$ with probability $p(A\cap[I])/p(A)$. This corresponds
to the biased coin-tosses in Lemma \ref{L2}. And the fair
coin-tosses in Lemma \ref{L2} reflect the random ordering of branches used in defining the
depth-first order, as can be seen from the definition in Section
\ref{con-ptrees} 
(the set of children of any vertex is put in exchangeable random order). The
only exception is on children of $j$ itself, which are all in ${\cal N}^*(v)$,
so the $b_l$'s are designated to artificially remove each of them with 
probability $1/2$. This establishes the lemma. \cq

The importance of the lemma is explained by the following formula
\begin{equation}\lb{asymG}
\max_v|\gip(e(v))-p({\cal N}^*(v))-c^*(v)|=o(\sp) \mbox{ in probability.}
\end{equation}
Since asymptotically we know that children of $i\in[I]$ are not in $[I]$,
and since by Lemmas \ref{limB} and \ref{limtheta}:
\begin{equation}\lb{negl}
\sp^{-1} \max_{v\notin[I]} p(\BB_v\setminus[I])
\cp 0,
\end{equation}
so in particular $\max_j\sp^{-1}c^*(j)\to 0$ in probability 
with the  notations
above, 
this is a straightforward consequence of Lemma \ref{limtheta} and
\begin{lmm}\label{keyg}
Suppose that no vertex $i\in [I]$ has a child that is also in $[I]$,
then we have for every $v$
\begin{equation}
\lb{keygip}
\gip(e(v))=p({\cal N}^*(v))-
\sum_{i\in{\cal N}(v)\cap[I]}(p_i-p({\cal B}_i)).
\end{equation}
\end{lmm}

\proof
Recall by definition \re{rhok} of the processes $\rho_k$ that if
$k$ is a child of some $i\in[I]$,
$\rho_k(e(v))=p_k$ whenever $v$ is examined after the parent $f(i)$
of $i$ {\em and}
strictly before $k$ in depth-first order, and $\rho_k(e(v))=0$ otherwise. 
As a consequence of \re{deprop}, we thus have
$$\gip(e(v))=p({\cal N}(v))-\sum_{i\in[I],k\in{\cal B}_i}p_k
\ind_{\{e(f(i))\leq e(v)<e(k)\}}.$$
A careful examination of this formula shows that a term in the sum on the
right is not zero if either $v$ has some ancestor $i\in[I]$, or some ancestor
of $v$ has a child $i\in [I]$ that is after $v$ in depth-first order,
and these situations are exclusive by the assumption that vertices of
$[I]$ do not have children in $[I]$.
In the first case, the formula says that we remove all the $\bp$-values
of children of $i$ that are after $v$ in depth-first order, in the
second case, it says that we remove the $\bp$-values of all the
children of $i$, implying \re{keygip}. \cq

{\em Proof of Proposition \ref{P1}.}
Fix $\eps > 0$ and consider arbitrary $v \in [n]$. Recall the
definition of ${\cal A}(v), {\cal G}(v), {\cal N}^*(v), c^*(v)$. 
We assert, from
Lemmas   \ref{L2}   and   \ref{L3},   that   for  any   $K>0$   there   exists
$\delta=\delta(\eps,K)$ with
\begin{equation}
P\left( |p({\cal N}^*(v))-c^*(v) - \sfrac{1}{2}{\cal G}(v) | > \eps \sp
| {\cal A}(v) \right)
\leq \exp(-\delta/\sp)\mbox{ on }
\{{\cal G}(v) \leq (K+2\eps)\sp\} . \lb{FGv}
\end{equation}
To argue (\ref{FGv}), note that conditioning
on the set $A={\cal A}(v)$ of vertices in the path from the root to $v$
determines the value
${\cal G}(v) := p({\cal A}(v)\setminus[I]) = \bar{q}_1$ say.
Then Lemmas \ref{L2}, \ref{L3} 
imply that the conditional distribution of $p({\cal N}^*(v))-c^*(v)$
has 
the
distribution of $X$ in Lemma \ref{L2}, 
The conclusion of Lemma \ref{L2} now gives (\ref{FGv}).

So for fixed $K$ and arbitrary $v \in [n]$
\[
P\left( |p({\cal N}^*(v))-c^*(v) - \sfrac{1}{2}{\cal G}(v)| > \eps \sp, \
\sfrac{1}{2}{\cal G}(v) \leq (K+2\eps) \sp \right)
\leq \exp(-\delta/\sp)= o(1/n) . 
\]
Using Boole's inequality gives 
$$P\left(\sfrac{1}{\sp}|p({\cal    N}^*(v))-c^*(v)-\sfrac{1}{2}{\cal    G}(v)|>\eps
\mbox{ for some }v \mbox{ with }\sfrac{1}{2\sp}{\cal G}(v)\leq K+2\eps\right)=o(1).$$
By (\ref{asymG}) we may replace $p({\cal N}^*(v))-c^*(v)$
by $\gip(e(v))$ in the previous expression.
We now use a slightly fussy truncation procedure. Imposing an extra constraint,
\begin{eqnarray}
&&
P\left(
\sfrac{1}{\sp} \max_v \gip(e(v)) \leq K,
\sfrac{1}{\sp}|\gip(e(v)) - \sfrac{1}{2}{\cal G}(v)| > \eps
\mbox{ for some } v \right. \nonumber\\
&& \left. \mbox{ with }
\sfrac{1}{2\sp} {\cal G}(v) \leq K+2\eps
\right) = o(1) .
\lb{FGK}
\end{eqnarray}
We claim that we can remove the restriction on $v$ to get
\begin{equation}
P\left(
\sfrac{1}{\sp} \max_v \gip(e(v)) \leq K,
\sfrac{1}{\sp}|\gip(e(v)) - \sfrac{1}{2}{\cal G}(v)| > \eps
\mbox{ for some } v \right)
= o(1) .
\lb{FGK2}
\end{equation}
Indeed, if $v$ has parent $v^\prime$ then
${\cal G}(v) - {\cal G}(v^\prime) \leq \max_{i\notin[I]} p_i =o(\sp)$. 
So if there exists a $v$ with
$\sfrac{1}{2\sp} {\cal G}(v) > K+2\eps$
then (for large $n$) there is an ancestor $w$ with
$K+\eps < \sfrac{1}{2\sp} {\cal G}(w) < K+2\eps$.
But if the first event in
(\ref{FGK}) occurs, one obviously cannot have
$\sp^{-1}|\gip(e(w))-\sfrac{1}{2}{\cal G}(w)|\leq \eps$ by definition of $w$.
Thus the probability in \re{FGK2} is bounded by twice the probability in
\re{FGK}.
This establishes (\ref{FGK2}).
Since Proposition \ref{pGZ} implies
$\sfrac{1}{\sp} \max_v \gip(e(v))$
is tight as $n \to \infty$,
(\ref{FGK2}) implies
\begin{equation}\label{cvgf}
\max_v \sfrac{1}{\sp} |\gip(e(v)) - \sfrac{1}{2}{\cal G}(v)| \cp 0 .
\end{equation}
Now let us show that
the sequence $(\sp\max_{v\in[n]}\HT(v),n\geq 1)$ is tight. Fix
$\eps>0$ and let
$K>0$ such that
$$P\left(\sfrac{1}{\sp}\max_v \gip(e(v))>K\right)<\eps/2,$$
Then
$$P\left(\sp\max_v\HT(v)>K+1\right)\leq \eps/2+P\left(\sp\max_v\HT(v)>K+1,
\sfrac{1}{\sp}\max_v \gip(e(v))<K\right),$$
but by the same kind of argument as above, if $\sp\max_v \HT(v)>K$, for $n$
large there must exist some $w$ with $K+1/2<\sp\HT(w)<K+1$. By Lemma \ref{L1}
we then have also $K+1/2<\sp^{-1}{\cal G}(v)<K+1$ with high probability,
so \re{cvgf} implies that the right-hand side in the last expression
is $<\eps/2$ for $n$ large.
This being proved, Lemma \ref{L1} rewrites as
$\max_v|\sp^{-1}{\cal G}(v)-\sp\theta_0^2\HT(v)|=o(1)$ in probability,
which together with (\ref{cvgf})
establishes the proposition. \cq

\section{Miscellaneous comments}

\noindent{\bf 1.} 
In principle Corollary \ref{httree} gives a criterion for boundedness
of $\TT^{\sth}$,
but one would prefer to have a condition directly in terms of $\btheta$.
Here are some steps in that direction.
From \cite[Theorem 1.1]{kal6715}, the process $\xbrt$
may be put in the form $\xbrt_s=X^1_s+X^2_s,s\geq 0$, where $X^1$ is a L\'evy
process on $[0,\infty)$ and $X^2$ has exchangeable increments
on $[0,1]$ and in a certain sense behaves less wildly than $X^1$. Precisely,
$X^1$ has no drift, its Gaussian part is $\theta_0$ and its L\'evy measure
is $\Lambda(\d x)=\sum_{i\geq 1}\delta_{\theta_i}(\d x)$, where $\delta_y(\d
x)$
is the Dirac mass at $y$. On the other hand, $X^2$ can be put in the form
$$X^2_s=-X^1_1 s + \sum_{i\geq 1}\tau_i(\ind_{\{s\leq V_i\}}-s)$$
for some square-summable random family $(\tau_i)$ and a sequence $V_i$
of independent r.v.'s with uniform law (notice
that $X^1$ and $X^2$ are by no means independent). Then, writing
$\kappa_{\xbrt}=\inf\{c>0:\sum_{i\geq 1}\theta_i^c<\infty\}$ and
$\kappa_{X^2}=\inf\{c>0:\sum_{i\geq 1}\tau_i^c<\infty\}$
we have that
\begin{equation}
\kappa_{X^2}\leq \frac{\kappa_{\xbrt}}{1+\sfrac{1}{2}\kappa_{\xbrt}},
\end{equation}
which is what we mean by ``behaving less wildly''.
It is therefore reasonable that the
problem on the finiteness of the integral
$\int^1 \d s/X^{\sth}_s$, which is a problem dealing with the behavior
at the left of the overall minimum of $\xbrt$, should be
replaced by a problem on the L\'evy process $X^1$ as soon as one can show
that the overall minimum of $\xbrt$ is actually attained at a local minimum
of $X^1$, and such that locally $X^2$ is negligible compared to $X^1$ at this
time. Since $X^1$ has no negative
jumps, the time-reversed process has no positive jumps, and such questions are
addressed in Bertoin \cite{bertoin95lrg} and Millar
\cite{millar81c}. Pushing the intuition one  step further, by analogy with the
standard criterion for  non-extinction of continuous-state branching processes
and the analogy of ICRT's and L\'evy trees mentioned above, we conjecture that
$\int^{\infty}\Psi^{-1}(\lambda)\d   \lambda<\infty$  is  equivalent   to  the
boundedness of $\TT^{\sth}$, where $\Psi$ is the Laplace exponent of $X^1$:
$$\Psi(\lambda)=\theta_0\lambda^2/2+\sum_{i\geq 1}
(\exp(-\lambda \theta_i)-1+\lambda\theta_i). $$

\noindent{\bf 2.} 
As we mentioned before, a natural  guess would be that the exploration process
of $\TT^{\sth}$ in the general case $\theta_0>0$ is
$\sfrac{2}{\theta_0^2}Y^{\sth}$. 
It is more difficult to get an intuition of what the exploration process of 
$\TT^{\sth}$ should be in the cases when $\theta_0=0$, when the Brownian part 
of $X^{\sth}$ vanishes. By the general theory of continuum random trees, it
should be easy to prove that compactness of the tree is enough to obtain
the existence of an exploration process for $\TT^{\sth}$, which is the
weak limit of $2(\theta^n_0)^{-2}Y^{\sth^n}$ for some
$\btheta^n\in \bTheta\to\btheta$ pointwise with
$\theta^n_0>0$ for every $n$. But this would not tell much about the look of
this process. Another way would be to try to generalize local time methods
used in \cite{DGall02}, but these do not seem to adapt so easily to bridges
with exchangeable increments instead of L\'evy processes.

\noindent{\bf Acknowledgment.  }Thanks to an anonymous referee for a careful
reading of a former version of the paper.


\begin{thebibliography}{99}


\bibitem{me46}
D.J.\ Aldous: The continuum random tree I. {\em Ann.\ Probab.\ }{\bf 19},
1--28 (1991).
\MR{1085326}

\bibitem{me55}
D.J.\ Aldous: The continuum random tree II: an overview. In M.T.\ Barlow and
N.H.\ Bingham, editors, Stochastic Analysis, pp.\ 23--70. Cambridge University
Press (1991).
\MR{1166406}

\bibitem{me56}
D.J.\ Aldous: The continuum random tree III. {\it
Ann.\ Probab.\ } {\bf 21},248--289 (1993).
\MR{1207226}

\bibitem{me102}
D.J. Aldous, G. Miermont and J. Pitman:
Brownian Bridge Asymptotics for Random $p$-Mappings.
{\it Electr. J. Probab.}
{\bf 9}, 37-56, (2004).
\article{http://www.math.washington.edu/~ejpecp/EjpVol9/paper3.abs.html}


\bibitem{meAMP}
D.J. Aldous, G. Miermont and J. Pitman:
Weak Convergence of Random $p$-Mappings and the Exploration
Process of the Inhomogeneous Continuum Random Tree.
\arxiv{math.PR/0401115}

\bibitem{me61}
D.J.\ Aldous and J.\ Pitman: Brownian bridge asymptotics for random mappings.
{\em Random Structures Algorithms} {\bf 5}, 487--512 (1994).
\MR{1293075}

\bibitem{me82}
D.J.\ Aldous and J.\ Pitman: The standard additive coalescent.
{\it Ann.\ Probab} {\bf 26}, 1703--1726 (1998)
\MR{1675063}

\bibitem{me88}
D.J.\ Aldous and J.\ Pitman: A family of random trees with random
edge-lengths. {\em Random Structures Algorithms} {\bf 15} 176--195 (1999).
\MR{1704343}

\bibitem{me87}D.J.\ Aldous and J.\ Pitman: Inhomogeneous continuum
random trees and the entrance boundary of the additive coalescent.
{\it Probab.\ Th.\ Rel.\ Fields} {\bf 118}, 455--482 (2000).
\MR{1808372}

\bibitem{me98}
D.J.\ Aldous and J.\ Pitman: Invariance principles for non-uniform random
mappings and trees. In V.\ Malyshev and A.M.\ Vershik, editors,
{\em Asymptotic Combinatorics with Applications in Mathematical Physics},
pp.\ 113--147. Kluwer Academic Publishers (2002).
\MR{1999358}

\bibitem{bertoin95lrg}
J.\ Bertoin: On the local rate of growth of L\'evy processes with no positive
jumps. {\em Stoch.\ Proc.\ App.\ } {\bf 55}, 91--100 (1995)
\MR{1312150}

\bibitem{bertoin00f}
J.\ Bertoin: A fragmentation process connected to Brownian
motion. {\it Probab.\ Theory Relat.\ Fields} {\bf 117}, 289--301 (2000).
\MR{1771665}

\bibitem{bertoin01eac}
J.\ Bertoin: Eternal additive coalescent and certain
bridges
with exchangeable increments. {\it Ann. Probab.\ } {\bf 29}, 344--360 (2001).
\MR{1825153}

\bibitem{by88e}
P.\ Biane and M.\ Yor: Sur la loi des temps locaux browniens pris en
un temps exponentiel. In {\em S\'eminaire de Probabilit\'es XXII}, pp.\
454--466. {\em Lecture Notes in Math.\ } {\bf 1321}, Springer (1988).
\MR{960541}

\bibitem{jpmc97b}
M.\ Camarri and J.\ Pitman: Limit distributions and random trees derived
from the birthday problem with unequal probabilities. {\em Electron.\
J.\ Probab.\ } {\bf 5}, Paper 2, 1--18 (2000).
\MR{1741774}

\bibitem{chasslou99}
P.\ Chassaing and G.\ Louchard: Phase transition for parking blocks, Brownian
excursion and coalescence. {\em Random Structures Algorithms} {\bf 21},
76--119 (2002).
\MR{1913079}

\bibitem{DGall02}
T.\ Duquesne and J.-F.\ Le Gall:
Random trees, L\'evy processes
and spatial branching processes. {\it Ast\'erisque} {\bf 281} (2002).
\MR{1954248}

\bibitem{durrett95}
R.\ Durrett. {\em Probability: theory and examples. } Duxbury Press, Belmont,
CA, second edition (1996).
\MR{1609153}

\bibitem{jeulin85}
T.\ Jeulin: Application de la th\'eorie du grossissement \`a
l'\'etude des temps locaux browniens. In T.\ Jeulin and M.\ Yor
(eds.): Grossissements de filtrations:
exemples et applications. {\em Lecture notes in Maths} {\bf 1118},
Springer, Berlin (1985).
\MR{884713}

\bibitem{kal73}
O.\ Kallenberg: Canonical representations and convergence criteria for
processes with interchangeable increments. {\em Z.\ Wahrsch.\ Verw.\ Gebiete}
{\bf 27}, 23--36 (1973).
\MR{394842}

\bibitem{kal6715}
O.\ Kallenberg: Path properties for
processes with independent and
interchangeable increments. {\em Z.\ Wahrsch.\ Verw.\ Gebiete}
{\bf 28}, 257--271 (1974).
\MR{402901}

\bibitem{kersting98m}
G.\ Kersting: On the height profile of a conditioned
Galton-Watson tree. 
\href{http://ismi.math.uni-frankfurt.de/kersting/research/profile.ps}{Preprint} (1998).

\bibitem{knight96e}
F.\ B.\ Knight: \rm The uniform law for exchangeable
and L\'evy process bridges. Hommage \`a P.\ A.\ Meyer et J.\
Neveu. {\it Ast\'erisque} {\bf 236}, 171--188 (1996).
\MR{1417982}

\bibitem{LeG93}
J.-F.\ Le Gall: The uniform tree in a Brownian excursion. {\it
Probab.\ Th.\ Rel.\ Fields} {\bf 96}, 369--383 (1993).
\MR{1231930}

\bibitem{LeG90}
J.-F.\ Le Gall: Brownian excursions, trees and measure-valued branching
processes. {\it
Ann.\ Probab.\ } {\bf 19}, 1299--1439 (1991).
\MR{1127710}

\bibitem{lglj97}
J.-F.\ Le Gall and Y.\ Le Jan: Branching processes in L\'evy processes:
The exploration process. {\em Ann.\ Probab.\ }{\bf 26}, 213--252 (1998).
\MR{1617047}

\bibitem{mm01}
J.-F.\ Marckert and A.\ Mokkadem: The depth-first processes of Galton-Watson
trees converge to the same Brownian excursion.
{\em Ann.\ Probab.} {\bf 31}, 1655--1678 (2003).
\MR{1989446}

\bibitem{miermont01}
G.\ Miermont: Ordered additive coalescent and fragmentations
associated to L\'evy processes with no positive jumps. {\it Electr.\ J.\
Probab.} {\bf 6}, paper 14 (33 pages) (2001).
\article{http://www.math.washington.edu/~ejpecp/EjpVol6/paper14.abs.html}
\MR{1844511}

\bibitem{millar81c}
P.W.\ Millar: Comparison theorems for sample function growth. {\em Ann.\
Probab.\ } {\bf 9}, 330--334 (1981).
\MR{606997}

\bibitem{op00}
C.A.\ O'Cinneide and A.V.\ Pokrovskii: Nonuniform random transformations.
{\em Ann.\ Appl.\ Probab.\ } {\bf 10} (4), 1151--1181 (2000).
\MR{1810869}

\bibitem{jp01hur}
J.\ Pitman: Random mappings, forests and subsets associated with
Abel-Cayley-Hurwitz multinomial expansions. {\em S\'eminaire Lotharingien de
Combinatoire}, Issue 46, 45 pp.\ (2001). 
\article{http://www.mat.univie.ac.at/~slc/wpapers/s46pitman.html}
\MR{1877634}

\bibitem{py92}
J.\ Pitman and M.\ Yor: Arcsine laws and interval partitions derived from a
stable subordinator. {\em Proc.\ London Math.\ Soc.\ } (3) {\bf 65}, 326--356
(1992).
\MR{1168191}

\bibitem{ry99}
D.\ Revuz and M.\ Yor: {\em Continuous martingales and
Brownian motion.} Springer, Berlin-Heidelberg, third edition (1999).
\MR{1725357}


\bibitem{ver79}
W.\ Vervaat: A relation between Brownian bridge and Brownian excursion.
{\em Ann. Probab.\ } {\bf 7}, 143--149 (1979).
\MR{515820}

\end{thebibliography}
\end{document}